\documentclass[11pt]{article}
\usepackage{amsmath}
\usepackage{amsthm}
\usepackage{bm}
\usepackage{multirow}
\newcommand{\supp}{\mathop{\mathrm{supp}}}
\usepackage{lipsum}
\usepackage{amsfonts}
\usepackage{graphicx}
\usepackage{epstopdf}
\usepackage{algorithm}
\usepackage{algorithmic}
\ifpdf
  \DeclareGraphicsExtensions{.eps,.pdf,.png,.jpg}
\else
  \DeclareGraphicsExtensions{.eps}
\fi

\theoremstyle{remark}

\theoremstyle{proposition}
\newtheorem{proposition}{Proposition}[section]

\newcommand{\TheTitle}{An efficient multigrid solver for isogeometric analysis}

\title{{\TheTitle}
\footnote{Francisco J. Gaspar has received funding from the European Union's Horizon 2020 research and innovation programme under the Marie Sklodowska-Curie grant agreement No 705402, POROSOS. The work of Carmen Rodrigo is supported in part by the Spanish project FEDER /MCYT MTM2016-75139-R and the Diputaci\'on General de Arag\'on (Grupo de referencia APEDIF, ref. E24\_17R).}}

  \author{
  \'Alvaro P\'e de la Riva\footnote{IUMA and Applied Mathematics Department, University of Zaragoza, Spain,(apedelariva@gmail.com)}
  \and
  Carmen Rodrigo\footnote{IUMA and Applied Mathematics Department, University of Zaragoza, Spain,(carmenr@unizar.es)}
  \and 
  Francisco J. Gaspar\footnote{CWI, Centrum Wiskunde and Informatica, Amsterdam, The Netherlands, (gaspar@cwi.nl)}
}

\usepackage{amsopn}
\DeclareMathOperator{\diag}{diag}

\newtheorem{definition}{Definition}

\begin{document}

\maketitle

\begin{abstract}
The design of fast solvers for isogeometric analysis is receiving a lot of attention due to the challenge that offers to find an algorithm with a robust convergence with respect to the spline degree. Here, we analyze the application of geometric multigrid methods to this type of discretizations, and we propose a multigrid approach based on overlapping multiplicative Schwarz methods as smoothers. The size of the blocks considered within these relaxation procedures is adapted to the spline degree. A simple multigrid V-cycle with only one step of pre-smoothing results to be a very efficient algorithm, whose convergence is independent on the spline degree and the spatial discretization parameter. Local Fourier analysis is shown to be very useful for the understanding of the problems encountered in the design of a robust multigrid method for IGA, and it is performed to support the good convergence properties of the proposed solver. In fact, an analysis for any spline degree and an arbitrary size of the blocks within the Schwarz smoother is presented for the one-dimensional case. The efficiency of the solver is also demonstrated through several numerical experiments, including a two-dimensional problem on a nontrivial computational domain.
\end{abstract}



\section{Introduction}\label{sec:1}
\setcounter{section}{1}

Isogeometric analysis (IGA) is a computational technique for the numerical solution of partial differential equations (PDEs), which was introduced by Tom Hughes et al. in the seminal paper \cite{Hughes_CMAME2005}. Since then, this approach has been widely applied within different frameworks, and a detailed presentation of IGA together with a number of engineering applications can be found in the book \cite{Hughes}. IGA is based on the idea of using spline-type functions which are used in the computer aided design (CAD) software for the parametrization of the computational domain in order to approximate the unknown solution of the PDE. B-splines  or nonuniform rational B-splines (NURBS) are the most commonly used functions. There are several issues that make this approach advantageous over classical finite element methods (FEM). First, it allows to represent exactly some geometries like conic sections, and also more complicated geometries are represented more accurately by this technique than by traditional polynomial based approaches. In addition, this precise description of the geometry is incorporated exactly at the coarsest grid level, making unnecessary further communication with the CAD system in order to do a mesh refinement procedure, which moreover does not modify the geometry. Another important advantage is the higher continuity, since IGA provides up to $C^{p-1}$ inter-element continuity, denoting $p$ the polynomial order (see \cite{COTTRELL2007, Beirao2014}). This corresponds to the so-called isogeometric $k-$method, which is one of the three refinement strategies for IGA proposed in \cite{Hughes_CMAME2005}, together with the $h-$refinement (reducing the mesh size by knot insertion) and $p-$refinement (order elevation, i.e. increase of the spline degree). The $k-$refinement is unique to IGA and its main advantage is that it maintains the maximum possible smoothness $C^{p-1}$ for the spline space of degree $p$. Due to its high performance, this is the most popular refinement strategy in the IGA community, and this is the one studied in this work.





From the computational point of view, the efficient solution of the linear systems arising from the discretization of a PDE problem is a crucial point for its numerical simulation. When a discretization with high spline degrees is considered, this issue is even more challenging since the condition number of stiffness matrices grows exponentially with the spline degree. The study of the computational efficiency for direct and iterative solvers was initiated in the papers \cite{Collier2012,Collier2013}, respectively, and recently the design of iterative solvers has attracted much attention in the isogeometric community. Many efforts have been devoted to develop preconditioners for this type of discretizations (see for example \cite{Buffa2013, Beirao2012, Beirao2013,Gahalaut2013,Kleiss2012,Sangalli2016}). 

For classical finite element, finite difference or finite volume discretizations, multigrid methods \cite{Bra77, TOS01, Hackb}  are well-known to be among the fastest solvers showing optimal computational cost and convergence behavior. Thus, it seems natural to try to extend these methods to IGA, and in fact, in the early IGA literature, multigrid solvers for finite element methods have been directly transferred to isogeometric discretizations with only minimal adaptations. However, a naive application of these multigrid methods to the isogeometric case results in an important deterioration of the convergence of the algorithms when the spline degree is increased. In particular, multigrid methods based on standard smoothers, like Gauss Seidel smoother, are not robust with respect to the spline degree (see e.g. \cite{Gahalaut2013b}).  It was observed in \cite{Donatelli2015} that the spectral radius of the multigrid iteration matrices based on standard smoothers tends to one exponentially as $p$ increases.
As it was pointed out in \cite{DonatelliGMCS17}, this bad behaviour is due to the presence of many small eigenvalues associated with high-frequency eigenvectors. This deterioration of the convergence of standard multigrid algorithms has motivated advances toward robust multigrid methods with respect to the spline degree. In \cite{Donatelli2015} a multigrid method was constructed based on a preconditioned Krylov smoother at the finest level and in \cite{Hofreiher2015} a mass matrix was proposed as smoother within a multigrid framework. For both methods, an increase in the number of smoothing steps was needed in order to obtain robustness with respect to the spline degree. Also, a multigrid smoother based on an additive subspace correction technique, applying a different smoother to each of the subspaces, is studied in \cite{Hofreither2017_b, Hofreither2017}. In this work, we aim to propose a robust and efficient geometric multigrid algorithm for isogeometric discretizations, based on a family of overlapping multiplicative Schwarz-type methods as smoothers. We will show that by choosing an appropriate Schwarz-type smoother we can efficiently remove those high-frequency components of the error associated with the small eigenvalues. This makes possible to obtain a very simple and efficient solver by using a multigrid $V-$cycle with only one smoothing step. 

It is well-known that many details are open for discussion and decision in the design of a multigrid method for a target problem, since the performance of multigrid algorithms strongly depends on the choice of their components. There are no rules to facilitate this challenging task, but the local Fourier analysis (LFA) appears as a very useful tool for the design of the algorithm. The local Fourier analysis or local mode analysis was introduced by Achi Brandt in \cite{Bra77, Bra94}, and since then, it has become the main quantitative analysis for the convergence of multigrid algorithms. This analysis is based on the Fourier transform theory, and a good introduction can be found in \cite{Stu_Tro, TOS01, Wess}, and in the LFA monograph \cite{Wie01}. In particular, the application of LFA to analyze the smoothing properties of multiplicative Schwarz-type smoothers is not standard, and a special treatment is needed, see \cite{Scott_Kees, LFA_overlapping}. For the one-dimensional case, we provide an analysis for any spline degree and an arbitrary size of the blocks in the smoother. In this way, we can choose the suitable size of the blocks to perform in the Schwarz smoother depending on the spline degree.

The rest of the paper is structured as follows. In Section \ref{sec:2} the considered model problem is described, together with the basics of IGA. Section \ref{sec:3} is devoted to introduce the proposed multigrid method, making special emphasis on the description of the class of multiplicative Schwarz methods used as smoothers. The local Fourier analysis technique considered for the study of the convergence of the multigrid algorithm is introduced in Section \ref{sec:4}. The basics of LFA are included in this section, as well as the non-standard analysis necessary to study the proposed smoothers. At the end of Section \ref{sec:4}, some results obtained from the LFA are presented to choose, for each spline degree, an appropriate Schwarz smoother. Section \ref{sec:5} deals with the numerical results. One- and two-dimensional tests are presented, including an example with a complex geometry to show the real power of the isogeometric solver. Finally, in Section \ref{sec:6} some conclusions are drawn. 

\section{Preliminaries}\label{sec:2}
\setcounter{section}{2}

We consider the Poisson problem in $d$ spatial dimensions on the domain $\Omega = (0,1)^d$ with homogeneous Dirichlet boundary conditions,
\begin{equation}\label{model_problem}
\begin{array}{ccc}
-\Delta u  &=& f, \ \  \mbox{in} \quad \Omega, \\ 
u &=& 0, \ \  \mbox{on} \quad \partial \Omega.
\end{array}
\end{equation}
The variational formulation of problem \eqref{model_problem} is given by: find $u \in H_0^1(\Omega)$ such that 
$$
a(u,v) = (f,v), \quad \forall v \in H_0^1(\Omega),
$$
where
$$
a(u,v) = \int_{\Omega} \nabla u \cdot \nabla v \, {\rm d} x, \quad\mbox{and}\quad 
(f,v) = \int_\Omega f v \, {\rm d} x.
$$
Given a finite dimensional approximation space $V_h \subset H_0^1(\Omega)$, the Galerkin approximation of the variational problem reads as follows: find $u_h \in V_h$ such that 
\begin{equation} \label{approx_variational}
a(u_h,v_h) = (f,v_h), \quad \forall v_h \in V_h.
\end{equation}
If we fix a basis $\{\varphi_1, \ldots, \varphi_{n_h}\}$ for $V_h$, ${\rm dim} V_h = n_h$, then the solution of problem \eqref{approx_variational} can be written as $u_h = \sum_{i=1}^{n_h} u_i \varphi_i$. The coefficient vector ${\mathbf u} = (u_1,\ldots, u_{n_h})$ can be computed by solving the linear system
$
A {\mathbf u} = {\mathbf b},
$
where $A$ is the stiffness matrix obtained from the bilinear form $a(\cdot,\cdot)$, i.e. $A = (a_{i,j}) = (a(\varphi_j,\varphi_i))_{i,j=1}^{n_h}$, and ${\mathbf b}  = (f,\varphi_i)_{i=1}^{n_h}$ is the right-hand side vector. In finite element methods, the chosen approximation space $V_h$ is usually a space of continuous piecewise polynomials, whereas in the isogeometric analysis $V_h$ is a space of functions with higher continuity (up to order $p-1$, being $p$ the polynomial degree of the B-spline).
Next, we present the fundamental ideas and notation on isogeometric analysis.
 \subsection{Knot vectors and basis functions}
Univariate B-spline basis functions are constructed from knot vectors. 
\begin{definition}
	A knot vector $\Xi = \lbrace \xi_1 , \xi_2 , \ldots , \xi_{n+p+1} \rbrace$ is a non-decreasing sequence of real numbers called knots, where $n$ is the number of basis functions and $p$ is the polynomial degree of the B-spline curve. Every knot $\xi_i$ has a knot index $i \in \lbrace 1, \ldots , n+p+1 \rbrace$ and every interval $[\xi_i, \xi_{i+1})$, bounded by a pair of consecutive knots, is called knot span.
\end{definition}
Moreover, these knots can be understood as coordinates in the parameter space. If all knots are uniformly spaced the knot vector is called {\it uniform}. Otherwise, it is called {\it non-uniform}. When the first and the last knots are repeated $p + 1$ times, the knot vector is said to be {\it open}. An important property of open knot vectors is that the resulting basis functions are interpolatory at the ends of the parametric interval, but are not, in general, interpolatory at interior knots. In this work, we always use open knot vectors. Fixed a knot vector,  the B-spline basis functions  $N_{i}^p(\xi)$ are defined recursively by the Cox-de-Boor formula (see \cite{deboor}), starting with $p=0$ (piecewise constants).

\begin{definition}(Basis functions)
Given a knot vector $\Xi = \lbrace \xi_1 , \xi_2 , \ldots , \xi_{n+p+1} \rbrace$, the $i$th piecewise constant B-spline basis function $N_{i}^0 : [0,1] \rightarrow {\mathbb R}$, $ i=1,\ldots, n+p$, is defined as	
\begin{equation}
\label{constant_spline}
	N_{i}^0(\xi)=\left\lbrace \begin{array}{ll}
	1 & \text{if } \xi_{i} \leq \xi < \xi_{i+1}, \\
	0 & \text{otherwise.} \\
	\end{array}\right.
\end{equation}	
For every pair $(k,i)$ such that $1\leq k \leq p$, $1 \leq i \leq n+p-k$, the basis functions $N_{i}^k : [0,1] \rightarrow {\mathbb R}$ are given recursively by the Cox-de-Boor formula:
\begin{equation}
\label{order_spline}
	N_{i}^k(\xi)= \displaystyle\frac{\xi - \xi_{i}}{\xi_{i+k} - \xi_{i}} N_{i}^{k-1}(\xi) + \displaystyle\frac{\xi_{i+k+1} - \xi}{\xi_{i+k+1} - \xi_{i+1}} N_{i+1}^{k-1}(\xi),
\end{equation}	
in which fractions of the forma $0/0$ are considered  as zero.	
\end{definition}

\begin{figure}[h]
	\centering{
		$\begin{array}{c}
		\includegraphics[width=0.75\textwidth]{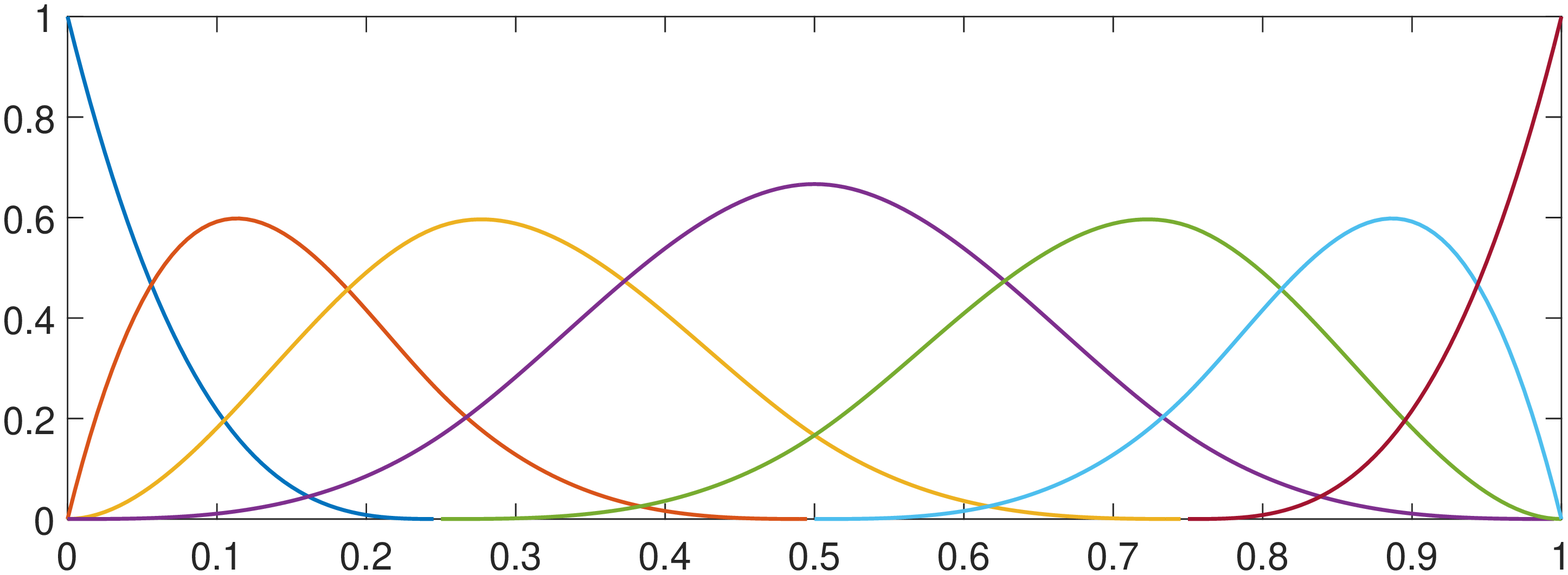}\\
		\end{array}$
	}
	\caption{Cubic B-splines basis functions with $\Xi = \lbrace 0,0,0,0,1/4,1/2,3/4,1,1,1,1 \rbrace$. }
	\label{cubic}
\end{figure}

\noindent B-spline basis functions own several important features and the reader can find them with their corresponding proofs in the literature \cite{Hughes,PiegTill96}. We remark some of them collected in the following proposition.
\begin{proposition}(B-spline basis functions properties)
	\begin{enumerate}
		\item For every $p$, the basis functions constitute a partition of unity
		$$\displaystyle\sum_{i=1}^{n} N_{i}^p(\xi) = 1, \forall \xi \in [\xi_1,\xi_{n+p+1}].$$ Clearly, out of the domain the basis functions are equal to $0$.
		\item B-spline basis functions are nonnegative over the domain and their support is always $p+1$ knot spans, $\supp(N_{i}^p) = [\xi_i,\xi_{i+p+1}]$.
		\item Across every knot $\xi_i$ the basis functions have $p-m_i$ continuous derivatives, with $m_i$ the multiplicity of knot $\xi_i$. 
	\end{enumerate} 
\end{proposition}
In this work we are interested in the isogeometric $k$-method \cite{Hughes}, which corresponds to spaces with maximum continuity, i.e. with $C^{p-1}$ regularity.  Thus, the knots considered in this work will have multiplicity one except for the first and  last knot. In Figure \ref{cubic} we illustrate a set of cubic B-splines basis for an open uniform knot vector.
Once the B-splines basis functions are defined, we introduce the notion of B-spline curve. 
\begin{definition}(B-spline curve)
	Given a knot vector $\Xi = \lbrace \xi_1 , \xi_2 , \ldots , \xi_{n+p+1} \rbrace$ we define a  B-spline curve  as the  linear combination	
	$${\mathbf C}(\xi) = \displaystyle\sum_{i=1}^{n} {\mathbf B}_i N_{i}^p(\xi),$$
	\noindent where ${\mathbf B}_i \in \mathbb{R}^2$ are the so-called control points of the curve, and $N_{i}^p$ are B-spline basis functions defined in \eqref{constant_spline}
and	\eqref{order_spline}.
\end{definition}
\subsection{Spline space}
As commented before, the solution of the variational formulation of problem \eqref{model_problem} is approximated in a spline space of degree $p$ with maximum smoothness. First, let us consider the simplest case when $d = 1$ (one dimensional problem). The computational domain then is the interval $\Omega = (0,1)$ and the corresponding two-point boundary value problem is:
$$
-u''(x) = f(x), \quad x \in \Omega, \quad  u(0) = u(1) = 0.
$$
Let the interval $(0,1)$ be subdivided into $m \in \mathbb{N}$ subintervals $I_i = ((i-1)h,ih)$, $i = 1, \ldots, m$, with $h=1/m$, and consider the knot sequence
$$
\Xi_{p,h} = \{\xi_1 = \ldots = \xi_{p+1} = 0 < \xi_{p+2} < \ldots \xi_{p+m} < 1 = \xi_{p+m+1} = \ldots
= \xi_{2p+m+1}\},
$$
where $\xi_{p+i+1} = i/m, i = 0, \ldots, m$. We define the spline space of degree $p \geq 1$ with maximum continuity
\begin{equation}
\label{spline_space}
S_{p,h}(0,1) = \{ u_h \in C^{p-1}(0,1) : u_h |_{I_i} \in {\mathbb P}^p, i = 1, \ldots, m, u_h(0) = u_h(1) = 0\}, 
\end{equation}
where $C^{p-1}(0,1)$ is the space of all $p-1$ times continuously differentiable functions on $(0,1)$, and ${\mathbb P}^p$ is the space of all polynomials of degree
less than or equal to $p$. The dimension of the space $S_{p,h}(0,1)$ is $p + m-2$, and the set of basis functions $\{N_{i}^p\}_{i=1}^{p+m-2}$ of this space are easily obtained by using \eqref{constant_spline}-\eqref{order_spline}. \\

In the case of higher spatial dimensions $d > 1$, the previous definitions are easily generalized  by means of tensor product. For simplicity, we assume that the spline degree $p$ and the number of subintervals $m$ are the same in all the directions.
However, this is not restrictive for the proposed solver described in this work.
In this way, over the domain $\Omega = (0,1)^d$ we define the spline space
$$
S_{p,h}(\Omega) = \underbrace{S_{p,h}(0,1) \oplus ... \oplus S_{p,h}(0,1)}_{d}.
$$
In the two-dimensional case, $d=2$, the knot vector 
$$
\Xi_{p,h} \times \Xi_{p,h} = \{(\xi,\eta), \xi \in \Xi_{p,h}, \eta \in \Xi_{p,h}\},
$$ generates a mesh of rectangular elements in the parametric space. The spline space is $S_{p,h}(0,1)^2 = S_{p,h}(0,1) \oplus S_{p,h}(0,1)$ and the bivariate B-splines basis is constructed by the tensor product of univariate B-splines basis. In this way a basis function $N_{i,j}^p: [0,1]^2 \rightarrow \mathbb{R}$ is defined in terms of
the univariate basis functions $N_i^p, N_j^p: [0,1] \rightarrow \mathbb{R}$ as follows:
$$
N_{i,j}^p(\xi,\eta) = (N_i^p \oplus N_j^p)(\xi,\eta) = N_i^p(\xi) N_j^p(\eta).
$$ 
In particular, for our model problem we can write 
$$
S_{p,h}(0,1)^2 = {\rm span} \{ N_{i,j}^p(\xi,\eta), i,j = 1, \ldots, p+m-2\}.
$$
\subsection{NURBS}
All the ideas presented here can be extended to solve partial differential equations in more complicated domains $\Omega$. In the isogeometric analysis, $\Omega$ is usually given by a NURBS parametrization. In this way, conic sections, such as circles and ellipses, can be represented exactly. NURBS are built through rational functions of B-splines. A NURBS basis function of degree $p$ is 
$$
R_i^p(\xi) = \frac{\omega_i N_i^p(\xi)}{\sum_{k=1}^{p+m} \omega_k N_k^p(\xi)},
$$
where $\{\omega_1, \ldots, \omega_{p+m}\}$ is a given set of weights. In this way,
a NURBS curve can be defined as the  linear combination
$${\mathbf C}(\xi) = \displaystyle\sum_{i=1}^{p+m} {\mathbf B}_i R_{i}^p(\xi),
$$
where ${\mathbf B}_i \in \mathbb{R}^2$ are the control points of the curve. In the two-dimensional case, NURBS basis functions are defined as
$$
R_{i,j}^p(\xi,\eta) = \frac{ \omega_{i,j} N_{i,j}^p(\xi,\eta)}{\sum_{k,l=1}^{p+m} \omega_{k,l} N_{k,l}^p(\xi,\eta) },
$$
and the domain $\Omega$ is represented by a geometry transformation ${\mathbf F}: (0,1)^2 \rightarrow \Omega$ given by
$$
{\mathbf F}(\xi,\eta) = \sum_{i=1}^{p+m} \sum_{j=1}^{p+m} {\mathbf B}_{i,j} R_{i,j}^p(\xi,\eta). 
$$
Now, the basis functions on $\Omega$ are defined by composing the basis functions on the parametric domain with the inverse of the geometry transformation, that is $R_{i,j} \circ {\mathbf F}^{-1}$,
and therefore the finite dimensional approximation space in our model problem is
$$
V_h = {\rm span} \{ R_{i,j} \circ {\mathbf F}^{-1}, i,j = 1, \ldots, p+m-2\}.
$$


\section{Multigrid method}\label{sec:3}
\setcounter{section}{3}


Multigrid methods are based on the smoothing property of a classical iterative algorithm and the acceleration of its convergence by a coarse-grid correction technique. The convergence of a multigrid method strongly depends on the choice of its components. Here, we consider the canonical prolongation operator and as the restriction its adjoint. Then, the coarse-grid operators will be constructed by Galerkin approximation. Regarding the type of cycle, we will demonstrate that $V-$cycles perform similarly as $W-$cycles, and therefore the former ones will be preferred. Moreover, we will also see that only one smoothing step will be enough to obtain a very efficient and robust multigrid for IGA. The choice of the smoother, however, will be more involved. 

\begin{figure}[htb]
\begin{center}
\begin{tabular}{cc}
\includegraphics[scale = 0.29]{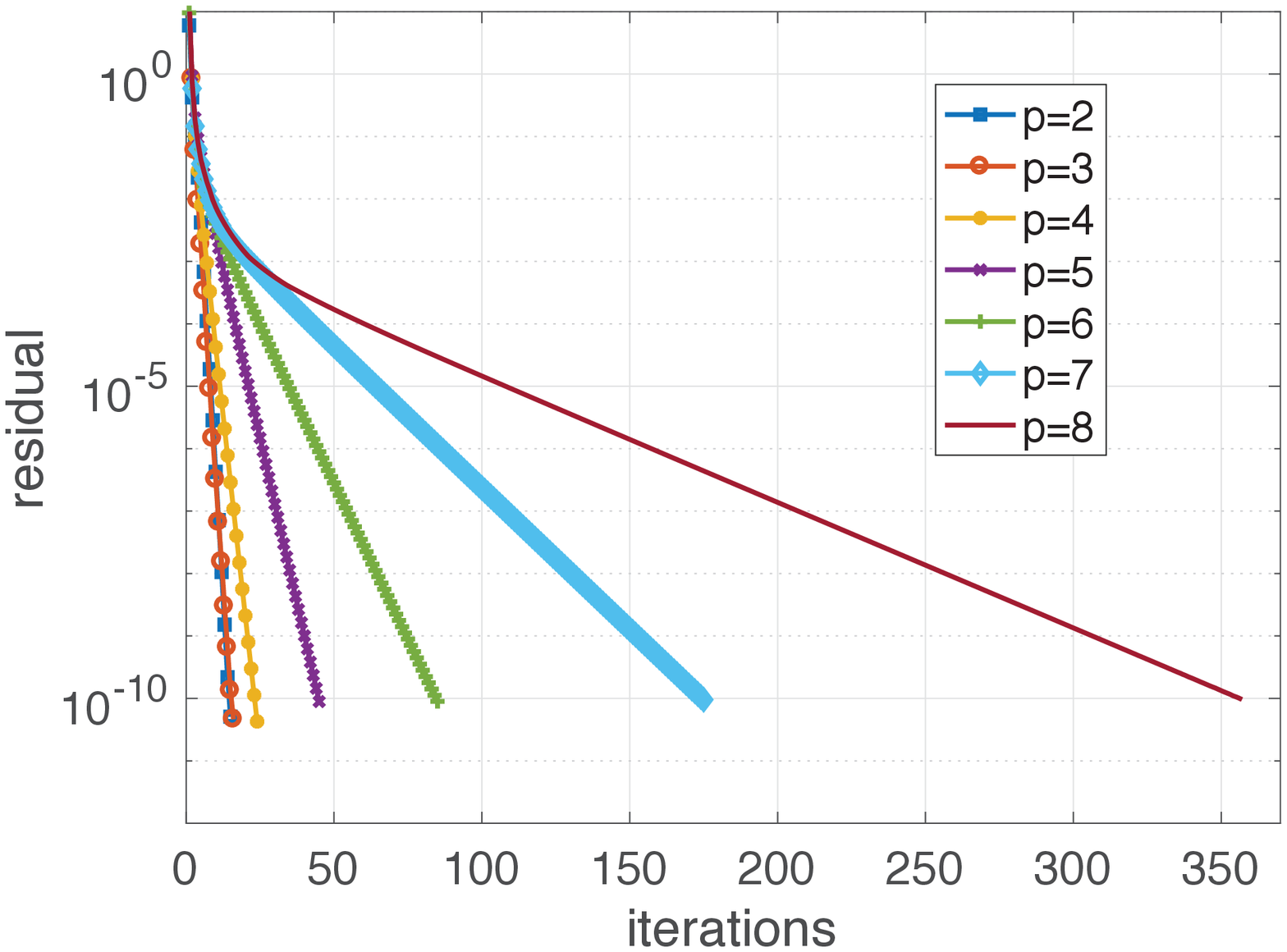}
&
\includegraphics[scale = 0.29]{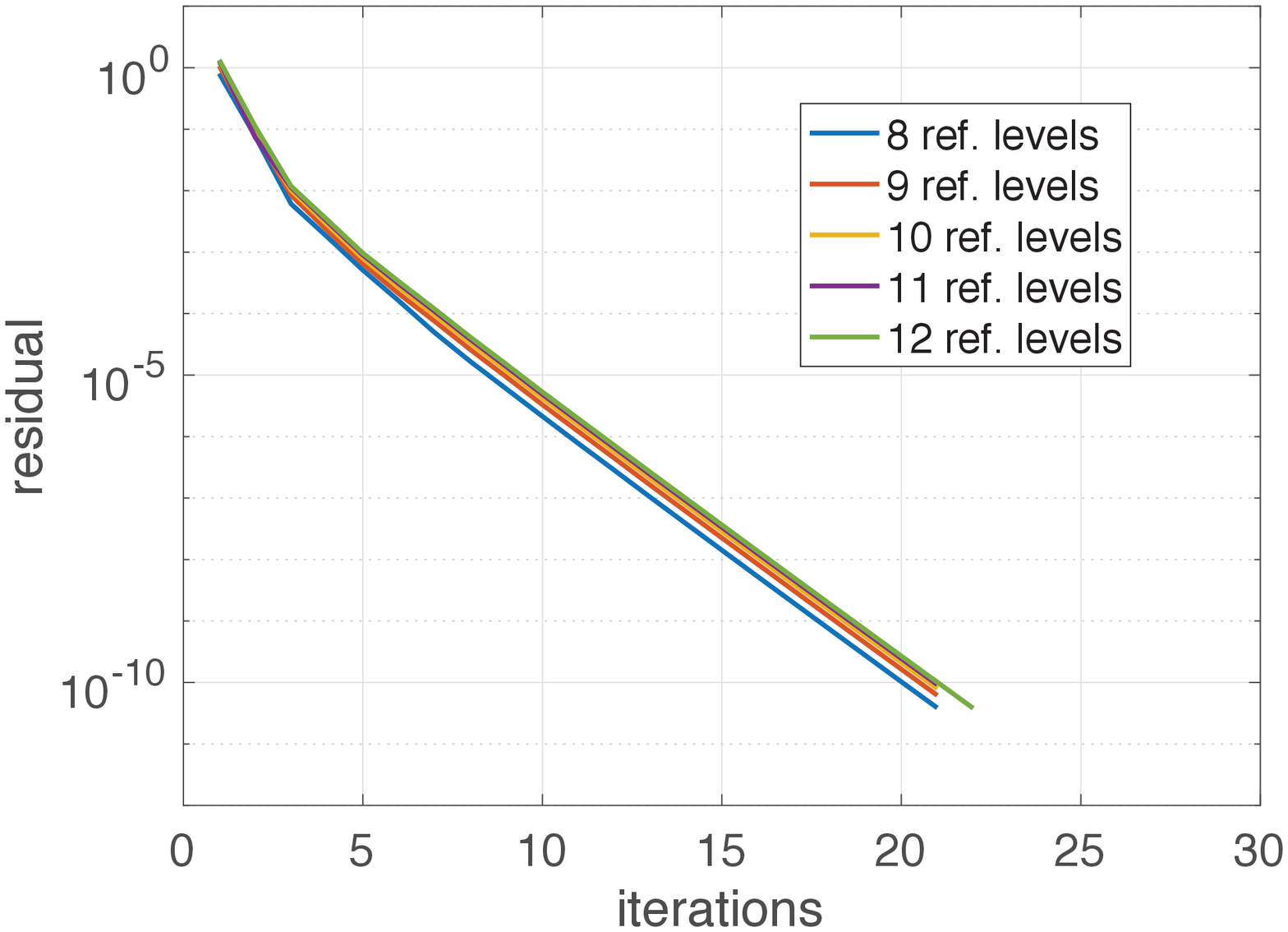}
\\
(a) & (b)
\end{tabular}
\label{MG_GS}
\caption{History of the convergence of a $V(1,0)-$multigrid method based on Gauss-Seidel relaxation for (a) different values of order $p$ and (b) for a fixed value $p=4$ and different grid sizes.}
\end{center}
\end{figure}

It is well-known that multigrid methods based on simple point-wise smoothers as Gauss-Seidel relaxation do not behave well when applied to isogeometric discretizations with larger values of $p$. This can be seen in Figure \ref{MG_GS} (a), where the number of iterations necessary to reduce the residual until $10^{-10}$ are displayed for different values of $p$ from $2$ to $8$. It is clear the deterioration of the multigrid convergence as soon as $p$ becomes larger. This implies that such a multigrid method is not reliable for large $p$, although its convergence is independent on the spatial discretization parameter for a fixed value of $p$. This latter can be seen in Figure   \ref{MG_GS} (b), where the history of the multigrid convergence is shown for $p=4$ with different target grids. 

Point-wise smoothers can be generalized to block-wise iterations by updating simultaneously a set of unknowns at each time, instead of only one. This is done by splitting the grid into blocks and solving together the equations corresponding to the grid-points in each block. There are many possibilities to construct these blocks. One can allow the blocks to overlap, giving rise to the class of overlapping block iterations, where smaller local problems are solved and combined via an additive or multiplicative Schwarz method.   

More specifically, the description of the multiplicative Schwarz iteration applied to the system $A{\mathbf u} = {\mathbf b}$ of order $n_h$ (see Section \ref{sec:2}) is as follows. Let us denote as $B$ the subset of unknowns involved in an arbitrary block of size $n$, that is $B = \left\{u_{k_1},\ldots, u_{k_n}\right\}$ where $k_i$ is the global index of the $i-$th unknown in the block. In order to construct the matrix to solve associated with such a block, that is $A^B$, we consider the projection operator from the vector of unknowns ${\mathbf u}$ to the vector of unknowns involved in the block. This results in a matrix $V_B$ of size $(n\times n_h)$, whose $i-$th row is the $k_i-$th row of the identity matrix of order $n_h$. Thus, matrix $A^B$ is obtained as $A^B = V_B A V_B^T$, and the iteration matrix of the multiplicative Schwarz method can be written as
$$\prod_{B=1}^{NB}\left( I-V_B^T(A^B)^{-1}V_B A \right),$$
where $NB$ denotes the number of blocks obtained from the splitting of the grid, which corresponds to the number of small systems that have to be solved in a relaxation step of the multiplicative Schwarz smoother. 

\begin{figure}[htb]
\begin{center}
\begin{tabular}{cc}
\includegraphics[width = 0.45\textwidth]{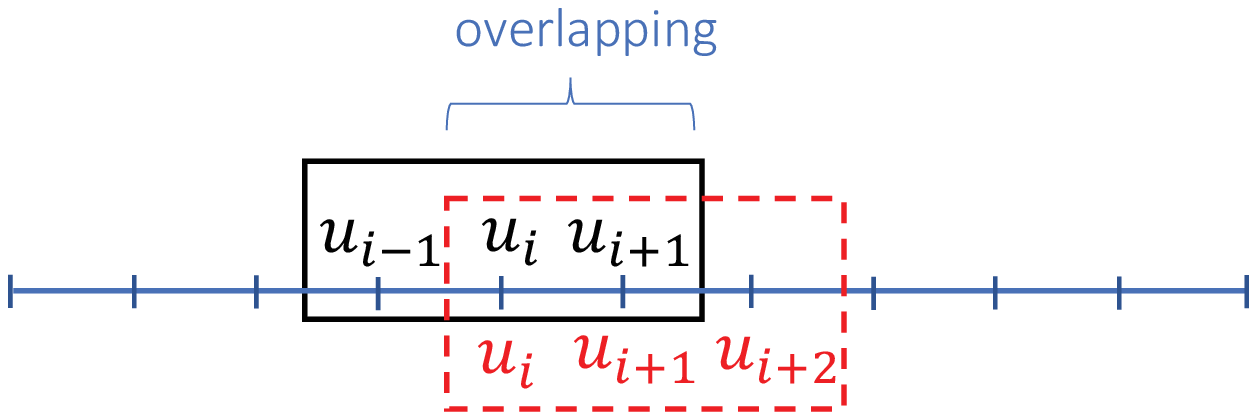}
&
\includegraphics[width = 0.45\textwidth]{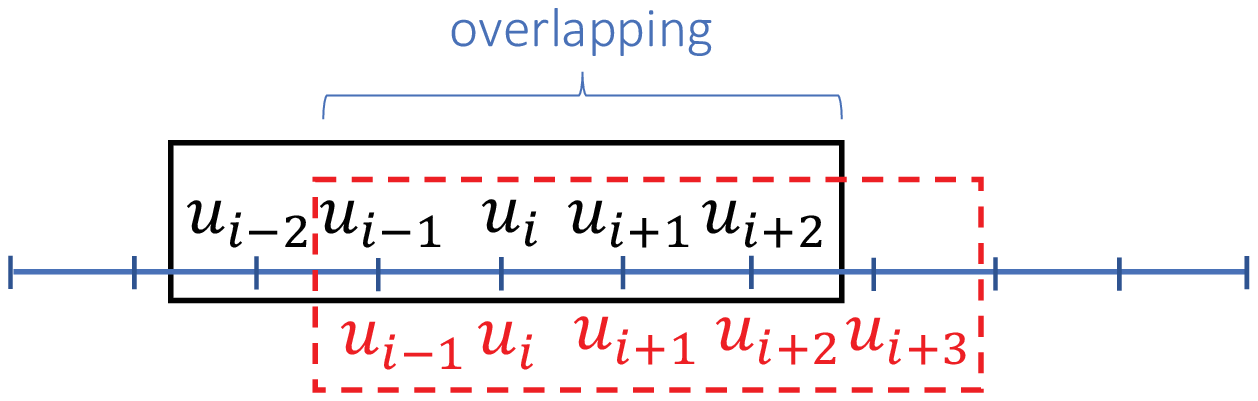}
\\ 
(a) & (b) 
\end{tabular}
\\
\vspace{-0.3cm}
\includegraphics[width = 0.5\textwidth]{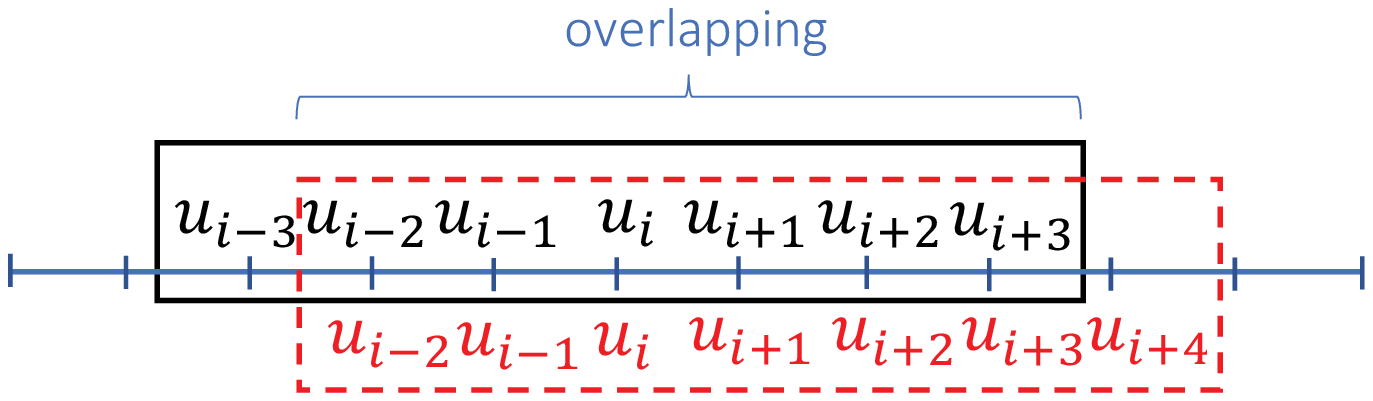}\\
(c)
\label{vankas_1D}
\caption{Size and overlapping of the blocks for the one dimensional smoothers: (a) Three- (b) Five- and (c) Seven-point multiplicative Schwarz iterations.}
\end{center}
\end{figure}

Here we consider multiplicative Schwarz methods with maximum overlapping.  In the one-dimensional case, our smoother will be based on blocks of three, five or seven points, depending on the spline degree $p$. Our study will be carried out up to $p=8$, but if one is interested in solving isogeometric discretizations with spline degree larger than $p=8$, only it is necessary to find the appropriate number of unknowns involved in the blocks to obtain an efficient multigrid approach. In the one dimensional case, the block of size $n$ associated with the $i-$th grid-point is composed of the unknown $u_i$, the $(n-1)/2$ unknowns on the right side of $u_i$ and the $(n-1)/2$ unknowns on its left side. More concretely, for the three-point multiplicative Schwarz smoother $B= \{u_{i-1}, u_i, u_{i+1}\}$, whereas for the seven-point relaxation, for example, we have $B= \{u_{i-3}, u_{i-2}, u_{i-1}, u_i, u_{i+1}, u_{i+2}, u_{i+3}\}$. In Figure \ref{vankas_1D} the corresponding blocks and the overlapping between them are schematized for the three-, five-, and seven-point multiplicative Schwarz smoothers (Figure \ref{vankas_1D} (a), (b) and (c), respectively). 
In the two dimensional case, we use the same idea, that is, we choose square blocks of size $\sqrt{n}\times \sqrt{n}$ around a grid-point $(i,j)$, yielding for example to the block $$B = \{u_{i-1,j-1}, u_{i,j-1}, u_{i+1,j-1}, u_{i-1,j}, u_{i,j}, u_{i+1,j}, u_{i-1,j+1}, u_{i,j+1}, u_{i+1,j+1} \}$$ for example for the nine-point multiplicative Schwarz smoother ($n=9$). In the same way, we will also consider the twenty five- and forty nine-point multiplicative Schwarz smoothers.   

In order to choose an appropriate multiplicative Schwarz smoother for each spline degree $p$, in the next section we introduce a local Fourier analysis that will help with this choice and will give a very good insight about why these relaxations work well.




\section{Local Fourier Analysis}\label{sec:4}  
\setcounter{section}{4}

In this section we introduce the basics of the local Fourier analysis and describe its application to IGA. First we will apply LFA to understand the difficulties in designing an appropriate multigrid for IGA and how this tool gives us insight to achieve this challenging task. After that, we will introduce LFA for overlapping multiplicative Schwarz smoothers. For simplicity in the presentation, we will describe the analysis for the one-dimensional problem. 

\subsection{Basics of LFA}
The main idea of local Fourier analysis is to assume a decomposition of the error function in Fourier modes, and then to study the behavior of each operator involved in the multigrid method on these components. Some assumptions should be considered to perform this analysis. In particular, LFA presumes that all the operations involved in the multigrid algorithm are local processes neglecting the effect of boundary conditions. 
Then, a regular infinite grid, ${\mathcal G}_h$, is assumed, which is obtained by the infinite extension of the considered spatial grid. 
By imposing some assumptions on the discrete operator $A_h$, as linearity and constant coefficients, a basis of complex exponential eigenfunctions of the operator, called Fourier modes or Fourier components, is obtained. They are defined as $\varphi_h(\theta,x) = e^{\imath \theta x}$, with $x\in{\mathcal G}_h$ and where $\theta\in\Theta_h:=(-\pi/h,\pi/h]$, yielding the so-called Fourier space ${\mathcal F}({\mathcal G}_h):=\hbox{span}\{\varphi_h(\theta,x)\,|\, \theta\in\Theta_h\}$. 
These Fourier components are divided into high- and low-frequency components on ${\mathcal G}_h$. We call low-frequency components to those Fourier modes associated with frequencies belonging to $\Theta_{2h} = (-\pi/2h,\pi/2h]$ (low frequencies), and the high-frequency components are those corresponding to the high frequencies $\theta\in \Theta_h\backslash \Theta_{2h}$. This classification depends on the coarsening strategy that we consider, which in this case is standard coarsening, that is, the coarse-grid  step size is obtained by doubling the step size of the fine grid. 

\begin{figure}[htb]
\begin{center}
\begin{tabular}{cc}
\includegraphics[scale = 0.33]{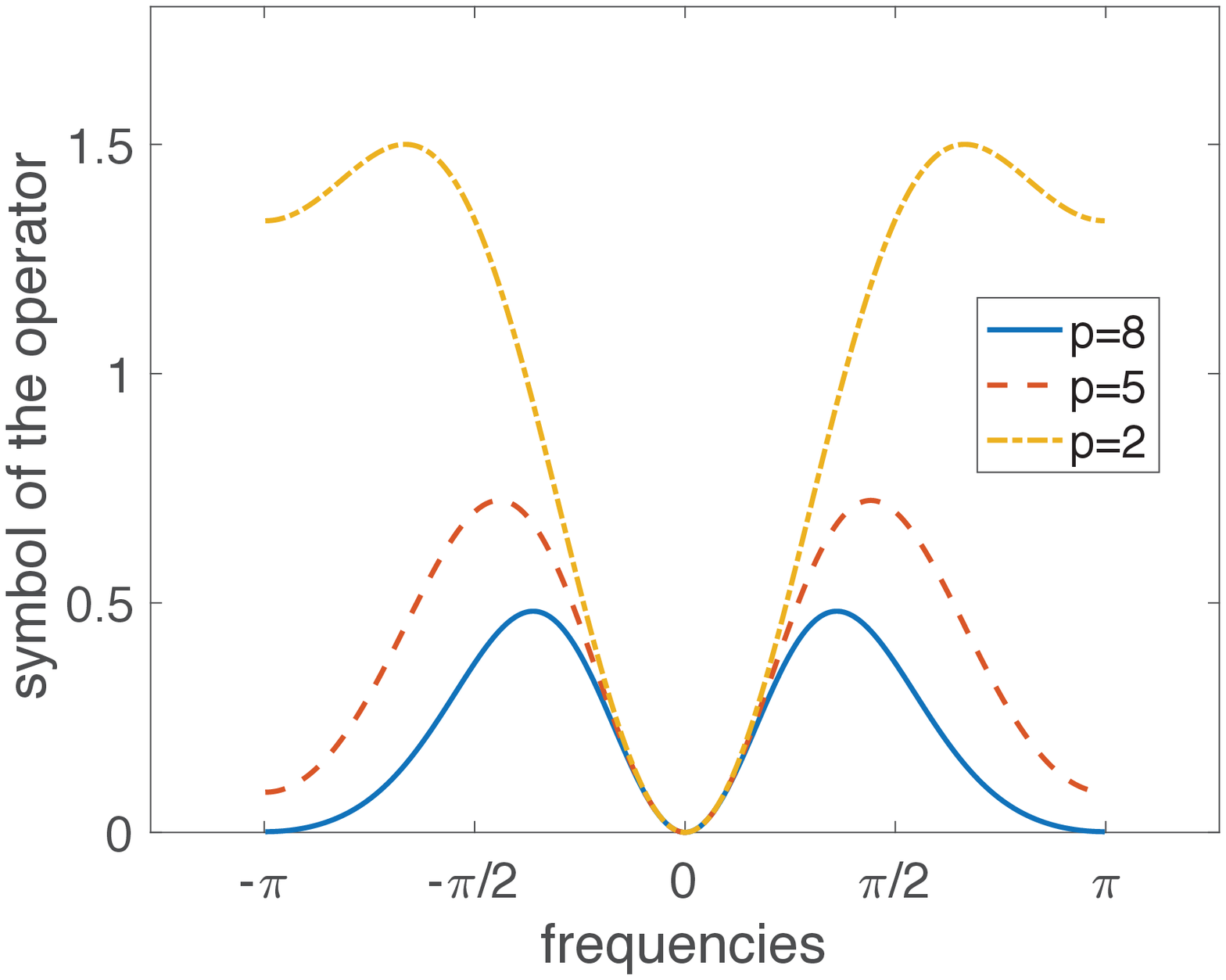}
&
\hspace{-0.4cm}
\includegraphics[scale = 0.33]{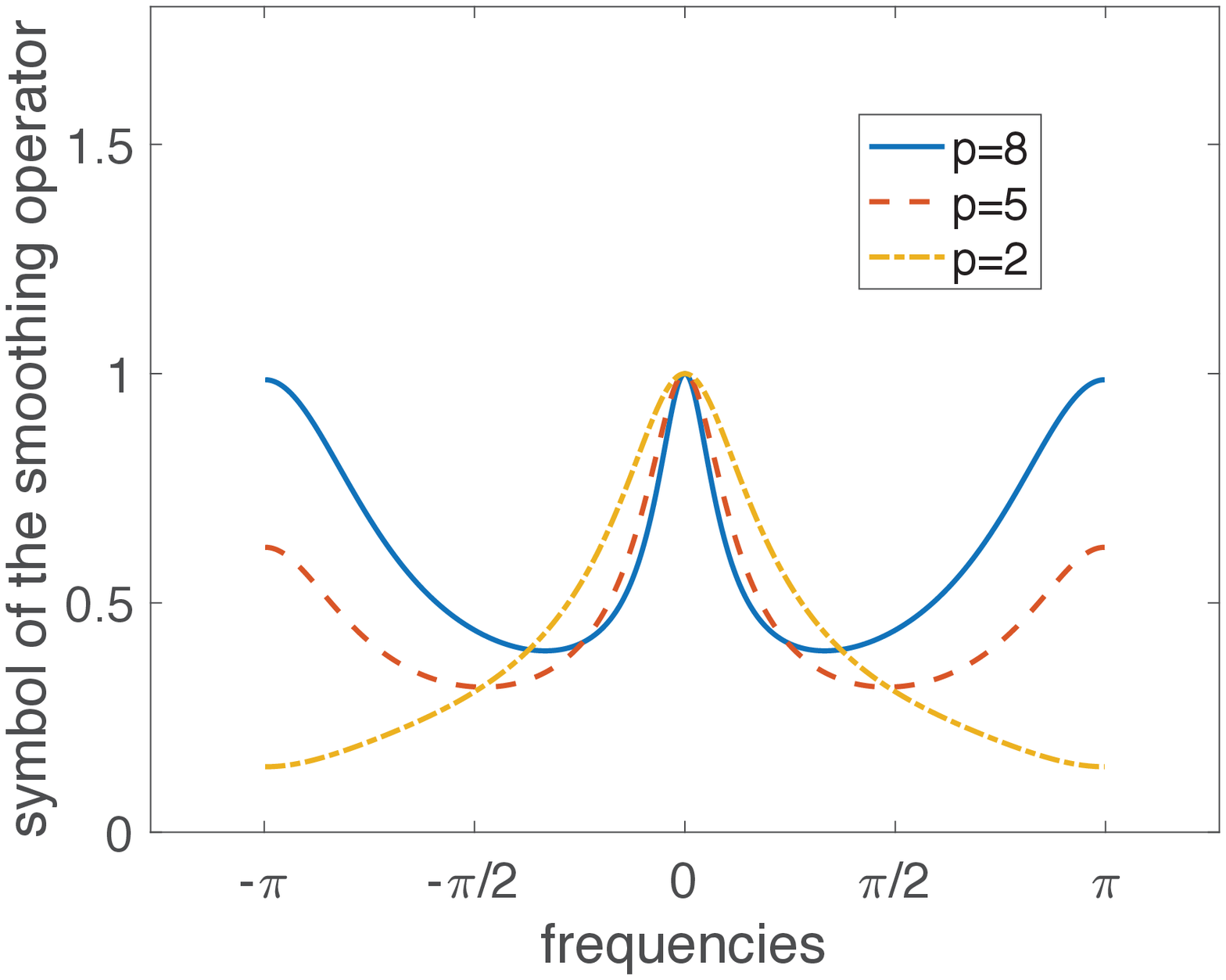}
\\
(a) & (b) 
\end{tabular}
\label{symbol_operator_GS}
\caption{Symbol of (a) the discrete operator and of (b) the Gauss-Seidel smoothing operator, for three different spline degrees $p = 2, 5, 8$.}
\end{center}
\end{figure}

Under the previous assumptions on the discrete operator, the Fourier components satisfy that $A_h \varphi_h(\theta,x) = \widetilde{A}_h(\theta)\varphi_h(\theta,x)$. This means that the Fourier components are ``eigenfunctions'' of the discrete operator, and the corresponding ``eigenvalues'' give rise to the so-called Fourier symbol of the operator, $\widetilde{A}_h(\theta)$. As an example,  the stencil corresponding to the one-dimensional IGA discrete operator with $p=2$ is, 
$$A_{h,2} = \frac{1}{h} \left[-\frac{1}{6}, -\frac{1}{3}, 1, -\frac{1}{3}, -\frac{1}{6} \right], $$
and its  Fourier symbol is given by
$$\widetilde{A}_{h,2}(\theta) = \frac{1}{h}\left(1-\frac{2\cos(\theta)}{3}-\frac{2\cos(2\theta)}{6}\right) = \frac{2}{3h}(2-\cos \theta (1+\cos \theta)).$$
In Figure \ref{symbol_operator_GS} (a) the eigenvalues of $A_{h,2}$ are displayed, together with those corresponding to the symbols of operators $A_{h,5}$ and $A_{h,8}$, that is, the discrete operators associated with $p=5$ and $p=8$, respectively. We can observe an important difference between the eigenvalues of $A_{h,2}$ and those for $A_{h,5}$ and $A_{h,8}$. Opposite to the case $p=2$, when $p=5$ or $p=8$, there are many small eigenvalues associated with high-frequency eigenvectors. This was also pointed out in  \cite{DonatelliGMCS17}, and it makes that the convergence of classical multigrid methods slows down since a standard smoother 
cannot provide an efficient damping of the components of the error associated with the high frequencies.  
In order to illustrate this, we analyze the Fourier representation of a Gauss-Seidel relaxation procedure. 
The Fourier components are also eigenfunctions of the smoothing iteration matrix $S_h$, which relates the error grid functions at two consecutive relaxation iterations. This means that we can obtain the error amplification factor or Fourier symbol of $S_h$ with respect to each frequency $\theta$, that is, $\widetilde{S}_h(\theta)$ (see \cite{TOS01, Wie01} for basic details). This is shown in Figure \ref{symbol_operator_GS} (b) for a lexicographic Gauss-Seidel smoother and the three different values of $p$ previously considered. It can be observed how this classical relaxation fails to annihilate the high-frequency components for $p=5$ and $p=8$, whereas it effectively reduces such Fourier modes for $p=2$. This is the reason why the careless application of a standard multigrid to IGA yields a deterioration of the convergence as the spline degree gets larger, as we previously reported in Figure \ref{MG_GS} (a).
The translation of the smoothing property, that is, the capacity of the smoother to eliminate the high-frequency components of the error, into a quantitative measure is the so-called smoothing factor, which is given as follows,
$$\mu = \sup_{\Theta_h\backslash \Theta_{2h}} \rho(\widetilde{S}_h(\theta)).$$
This measure, however, is not enough to analyze the interplay between the smoother and the coarse-grid correction technique that accelerates its convergence. In order to get more insight in the behavior of the multigrid method, taking into account the influence of the inter-grid transfer operators and the other components involved in the coarse-grid correction, at least a two-grid analysis is needed. The two-grid analysis consists of estimating the spectral radius of the two-grid operator. This error propagation operator is given as
$$M_h^{2h} = S_h^{\nu_2}(I_h-I_{2h}^hA_{2h}^{-1}I_h^{2h}A_h)S_h^{\nu_1},$$
where $A_h$ and $A_{2h}$ represent the discrete operators in the fine and coarse grids, $I_h^{2h}$ and $I_{2h}^h$ denote the restriction and prolongation operators, respectively, and $\nu_1$ and $\nu_2$ are the numbers of pre- and post-smoothing steps. It is well-known that the inter-grid transfer operators, and as a consequence the two-grid operator, couple some Fourier components. In particular, in the transition between fine and coarse grids, each low-frequency $\theta^{0}\in \Theta_{2h}$ is coupled with a high frequency $\theta^1 = \theta^0 - \hbox{sign}(\theta^0)\pi$, giving rise to the so-called spaces of $2h-$harmonics: ${\mathcal F}^2(\theta^{0}) = \hbox{span}\left\{\varphi_h(\theta^0,x), \, \varphi_h(\theta^1,x)\right\}$. This implies that the Fourier representation of the two-grid operator $M_h^{2h}$ with respect to ${\mathcal F}^2(\theta^{0})$, denoted by $\widetilde{M}_h^{2h}(\theta^0)$, is a $2\times 2$-matrix. We determine the spectral radius of $M_h^{2h}$ by calculating the spectral radius of these smaller matrices, that is, 
$$\rho_{2g} = \rho(M_h^{2h}) = \sup_{\theta^0\in\Theta_{2h}}\rho(\widetilde{M}_h^{2h}(\theta^0)).$$

Although the two-grid analysis is the basis for the classical asymptotic multigrid convergence estimates, a three-grid (or even a $k-$grid) analysis provides a deeper insight into the performande of multigrid. This analysis is crucial, for example, to study the behavior of $V-$cycles and the difference between the choice of pre- and post-smoothing steps. The error propagation matrix of a three-grid cycle can be written as follows:
$$M_h^{4h} = S_h^{\nu_2}(I_{h}-I_{2h}^h(I_{2h}-(M_{2h}^{4h})^{\gamma})A_{2h}^{-1}I_h^{2h}A_h)S_h^{\nu_1},$$
where $\gamma$ is the number of two-grid iterations, and $M_{2h}^{4h}$ is the two-grid operator between the two coarse grids, that is,
$$M_{2h}^{4h}=S_{2h}^{\nu_2}(I_{2h}-I_{4h}^{2h} A_{4h}^{-1}I_{2h}^{4h}A_{2h})S_{2h}^{\nu_1}.$$
In order to analyze how this three-grid operator acts on the Fourier modes, we take into account that not only in the transition from the finest to the second grid but also in the transition from the second to the coarsest grid there are some Fourier modes that are coupled. More concretely, four frequencies are coupled, and then we can define the so-called subspaces of $4h-$harmonics as
${\mathcal F}^4(\theta^0) = \hbox{span}\left\{\varphi_h(\theta^{\alpha}_{\beta}) \, | \, \alpha, \, \beta \in \{0,1\}\right\},$
where $\theta^0\in\Theta_{4h} = (-\pi/4h,\pi/4h]$ and
$\theta^{\alpha}_{\beta} = \theta^0 -\alpha\,\hbox{sign}(\theta^0)\pi/2 + (-1)^{\alpha+\beta}\beta\, \hbox{sign}(\theta^0)\pi.$ 
Thus, based on the decomposition of the Fourier space in terms of the subspaces of $4h-$harmonics, we can reduce the computation of the spectral radius of the three-grid operator to calculate the supreme of the spectral radii of the $4\times 4$ Fourier representations on these subspaces, $\widetilde{M}_h^{4h}(\theta^0)$, that is, 
$$\rho_{3g} = \rho(M_h^{4h}) = \sup_{\theta^0\in \Theta_{4h}}\rho(\widetilde{M}_h^{4h}(\theta^0)).$$

\begin{figure}[htb]
\begin{center}
\begin{tabular}{ccc}
\hspace{-0.5cm}
\includegraphics[scale = 0.21]{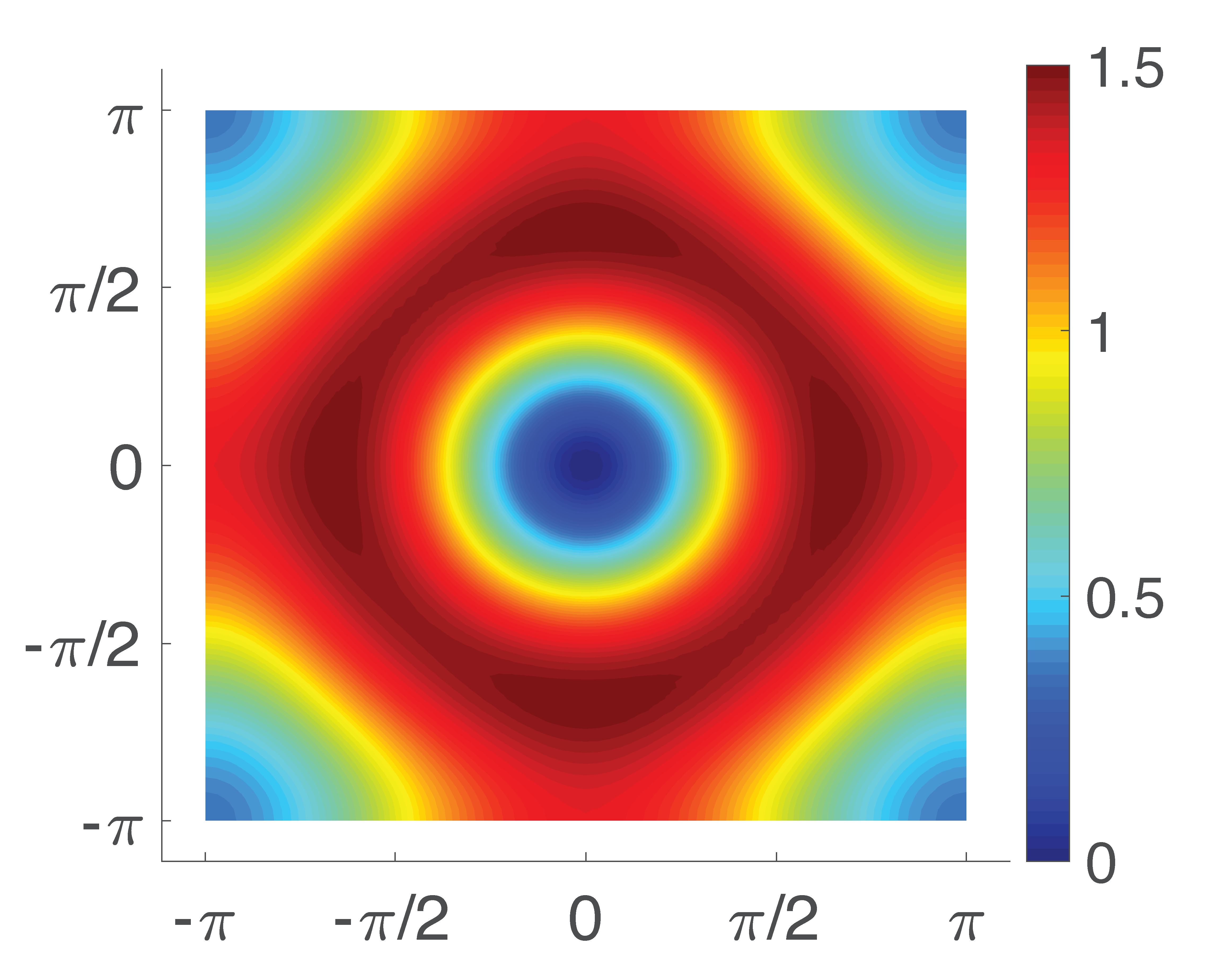}
&
\hspace{-0.5cm}
\includegraphics[scale = 0.21]{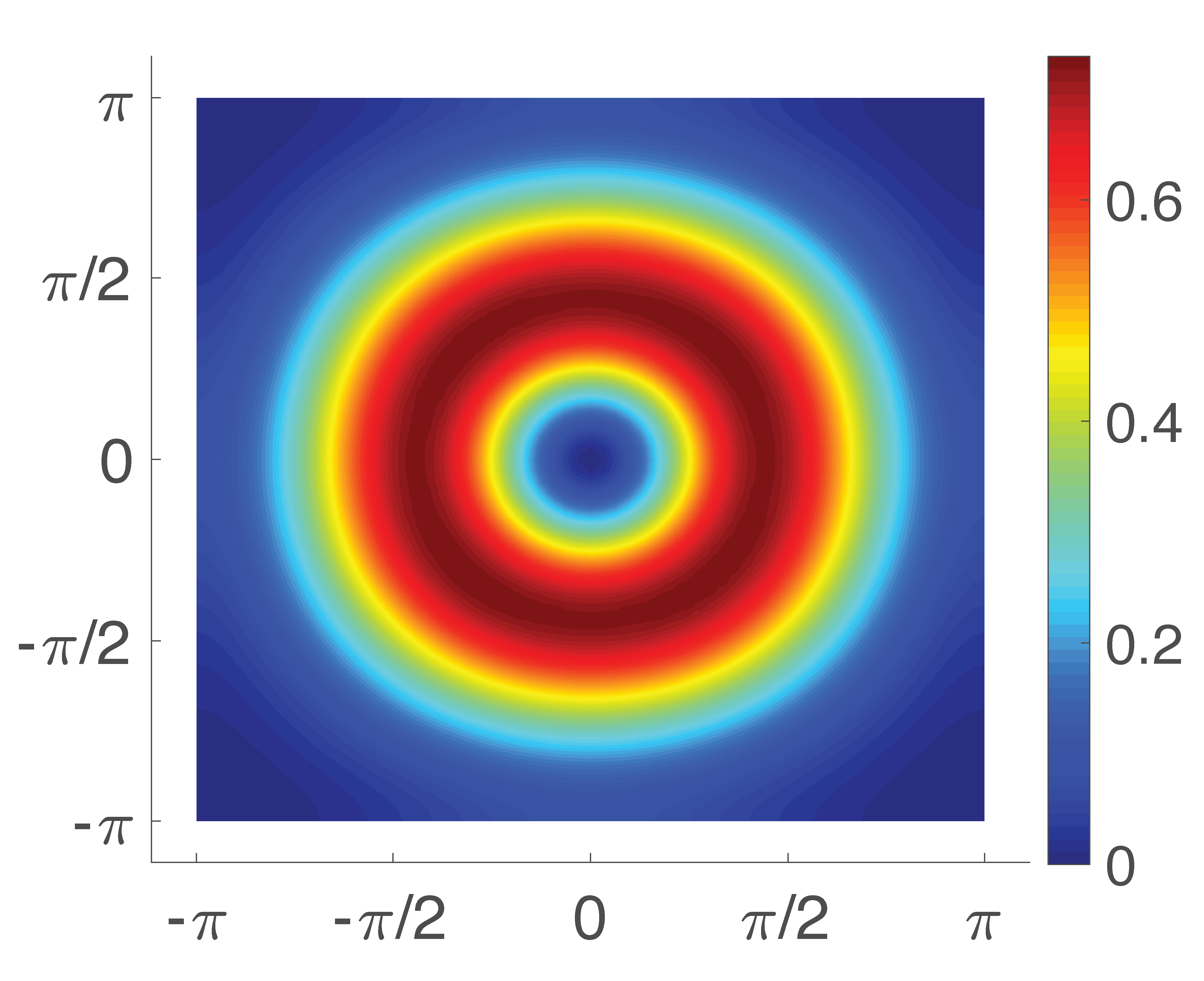}
&
\hspace{-0.5cm}
\includegraphics[scale = 0.21]{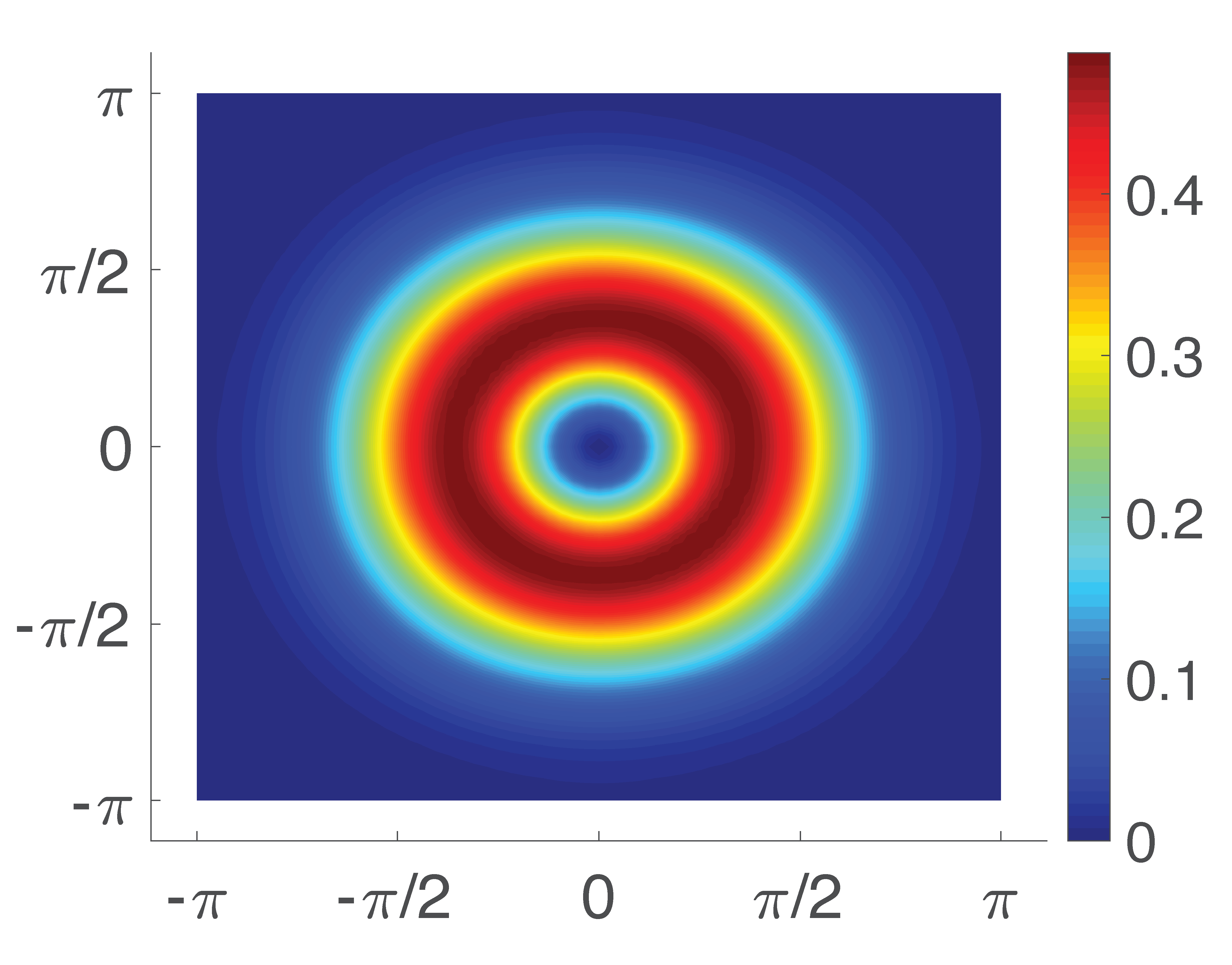}
\\
(a) & (b) & (c)
\end{tabular}
\label{symbol_operator_2D}
\caption{Symbol of the two-dimensional discrete operator for three different spline degrees (a) $p = 2$, (b) $p= 5$, and (c) $p= 8$.}
\end{center}
\end{figure}

This local Fourier analysis can be applied in two dimensions as well. In particular, we study the Fourier symbol of the discrete operator as we did before for the one-dimensional case. The results for three different spline degrees $p = 2, 5, 8$ are shown in Figure \ref{symbol_operator_2D}, where again we can observe a similar behavior to that obtained in 1D. For $p=2$ the eigenvalues associated with the high frequencies are mainly large, whereas as $p$ grows up there appear more and more small eigenvalues corresponding to the high frequency components. This implies that a standard two-dimensional smoother won't be able to completely annihilate the high-frequency components of the error and consequently a simple multigrid method will deteriorate its convergence as $p$ gets larger. 

\subsection{Local Fourier analysis for overlapping smoothers}\label{sec:4_schwarz}

In this section we describe the local Fourier analysis for the class of overlapping multiplicative Schwarz smoothers considered in this work. With this purpose, we follow the methodology described in \cite{LFA_overlapping} and the notation therein. This type of multiplicative Schwarz smoothers satisfy the invariance property, in the sense that the Fourier modes are their eigenvectors (see \cite{Scott_Kees} for a rigorous proof). The trouble, however, comes from their overlapping character. Due to this, the unknowns are updated more than once, what implies that, in addition to the initial and final errors, some intermediate errors appear, which have to be taken into account in the analysis. This makes necessary a special strategy to carry out the local Fourier analysis for overlapping block smoothers.
For simplicity in the presentation, we describe the analysis for the overlapping three-point multiplicative Schwarz smoother in 1D applied to the isogeometric discretization with $p=2$. At the end of this section we provide an expression to compute the symbol of the smoothing operator for any one-dimensional $n-$point multiplicative Schwarz relaxation for an arbitrary spline degree $p$.

For the considered case, at any grid-point $\xi_i$, three equations corresponding to the unknowns $u_{i-1}$, $u_i$, and $u_{i+1}$ are solved simultaneously. The block to solve, in terms of corrections and residuals, is given as follows,
\begin{equation}\label{system_block}
\frac{1}{h} \left( 
\begin{array}{ccc}
1 & -1/3 & -1/6 \\
-1/3 & 1 & -1/3 \\
-1/6 & -1/3 & 1
\end{array}
\right)
\left(
\begin{array}{c}
\delta u_{i-1} \\
\delta u_i \\
\delta u_{i+1} 
\end{array}
\right) = 
\left(
\begin{array}{c}
r_{i-1} \\
r_i \\
r_{i+1} 
\end{array}
\right).
\end{equation}
The corrections can be written in terms of the errors:
\begin{eqnarray*}
\delta u_{i-1} = e_h^{k+1}(\xi_{i-1}) - e_h^{k+2/3}(\xi_{i-1}),\\
\delta u_{i} = e_h^{k+2/3}(\xi_{i}) - e_h^{k+1/3}(\xi_{i}),\\
\delta u_{i+1} = e_h^{k+1/3}(\xi_{i+1}) - e_h^{k}(\xi_{i+1}),
\end{eqnarray*}
where $e_h^k$ denotes the initial error at $k-$iteration, $e_h^{k+1}$ represents the final error, and we define as $e_h^{k+1/3}$ and $e_h^{k+2/3}$ the intermediate errors appearing after the unknown has been updated once or twice, respectively. Thus, when solving the $i-$th block, for example, unknown $u_{i-1}$ has been already relaxed twice, and it will be updated the final third time. Without loss of generality, we consider that the error is given as a single Fourier mode multiplied by a coefficient $\alpha_{\theta}^{(m)}$, where $m=0,1,2,3$ is the number of times that the unknown has been updated in the current iteration. In this way, we can rewrite system \eqref{system_block} as follows,
\begin{eqnarray*}
&&\frac{1}{h} \left( 
\begin{array}{ccc}
1 & -1/3 & -1/6 \\
-1/3 & 1 & -1/3 \\
-1/6 & -1/3 & 1
\end{array}
\right)
\left(
\begin{array}{c}
\left( \alpha^{(3)}_{\theta}-\alpha^{(2)}_{\theta}\right)e^{-\imath \theta} \\
\left( \alpha^{(2)}_{\theta}-\alpha^{(1)}_{\theta}\right) \\
\left( \alpha^{(1)}_{\theta}-\alpha^{(0)}_{\theta}\right)e^{\imath \theta}
\end{array}
\right)\\ [1.3ex] 
&& \quad \quad=
-\frac{1}{h}\left(
\begin{array}{c}
-\frac{1}{6}\alpha^{(3)}_{\theta}e^{-\imath 3\theta}-\frac{1}{3}\alpha^{(3)}_{\theta}e^{-\imath 2\theta}+\alpha^{(2)}_{\theta}e^{-\imath \theta}-\frac{1}{3}\alpha^{(1)}_{\theta}-\frac{1}{6}\alpha^{(0)}_{\theta}e^{\imath \theta}\\ [1ex]
-\frac{1}{6}\alpha^{(3)}_{\theta}e^{-\imath 2\theta}-\frac{1}{3}\alpha^{(2)}_{\theta}e^{-\imath \theta}+\alpha^{(1)}_{\theta}-\frac{1}{3}\alpha^{(0)}_{\theta}e^{\imath \theta}-\frac{1}{6}\alpha^{(0)}_{\theta}e^{\imath 2\theta}\\ [1ex]
-\frac{1}{6}\alpha^{(2)}_{\theta}e^{-\imath \theta}-\frac{1}{3}\alpha^{(1)}_{\theta}+\alpha^{(0)}_{\theta}e^{\imath \theta}-\frac{1}{3}\alpha^{(0)}_{\theta}e^{\imath 2\theta}-\frac{1}{6}\alpha^{(0)}_{\theta}e^{\imath 3\theta}\\
\end{array}
\right).
\end{eqnarray*}
Since our aim is to find the relation between the initial and the fully corrected errors, we rearrange the previous system into a system of equations for the updated coefficients, that is,
\begin{equation*}
\underbrace{
\left(
\begin{array}{ccc}
-\frac{1}{6}e^{\imath \theta} & -\frac{1}{3} & e^{-\imath \theta}-\frac{1}{3}e^{-\imath 2\theta} -\frac{1}{6}e^{-\imath 3\theta}\\ [1ex]
-\frac{1}{3}e^{\imath \theta} & 1 & -\frac{1}{3}e^{-\imath \theta} -\frac{1}{6}e^{-\imath 2\theta}\\ [1ex]
e^{\imath \theta} & -\frac{1}{3} & -\frac{1}{6}e^{-\imath \theta}\\
\end{array}
\right)}_{P}
\left(
\begin{array}{c}
\alpha^{(1)}_{\theta} \\
\alpha^{(2)}_{\theta} \\
\alpha^{(3)}_{\theta}
\end{array}
\right)
= \underbrace{\left(
\begin{array}{c}
0 \\ [1ex]
\frac{1}{6}e^{\imath 2\theta} \\ [1ex]
\frac{1}{3}e^{\imath 2\theta}+\frac{1}{6}e^{\imath 3\theta}
\end{array}
\right)}_{Q}
\alpha^{(0)}_{\theta},
\end{equation*}
The amplification factor for the error is given by the last component of $P^{-1}Q$. Once that we have the symbol of the smoothing operator $\widetilde{S}_h(\theta) = (P^{-1}Q)_3$, the smoothing and $k-$grid local Fourier analysis can be carried out as explained before in the standard way.

\begin{figure}[htb]
\begin{center}
\begin{tabular}{cc}
\includegraphics[scale = 0.29]{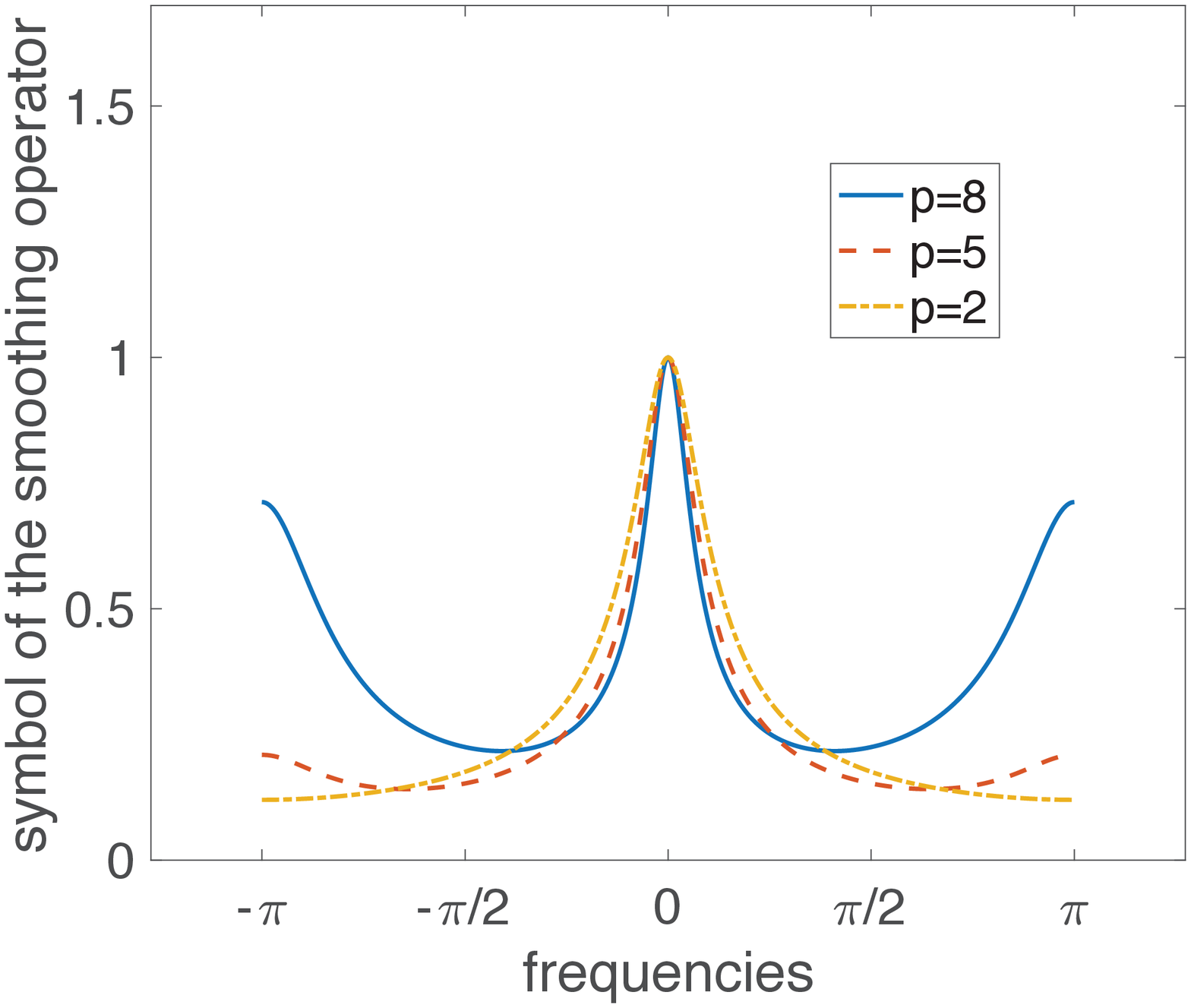}
&
\includegraphics[scale = 0.29]{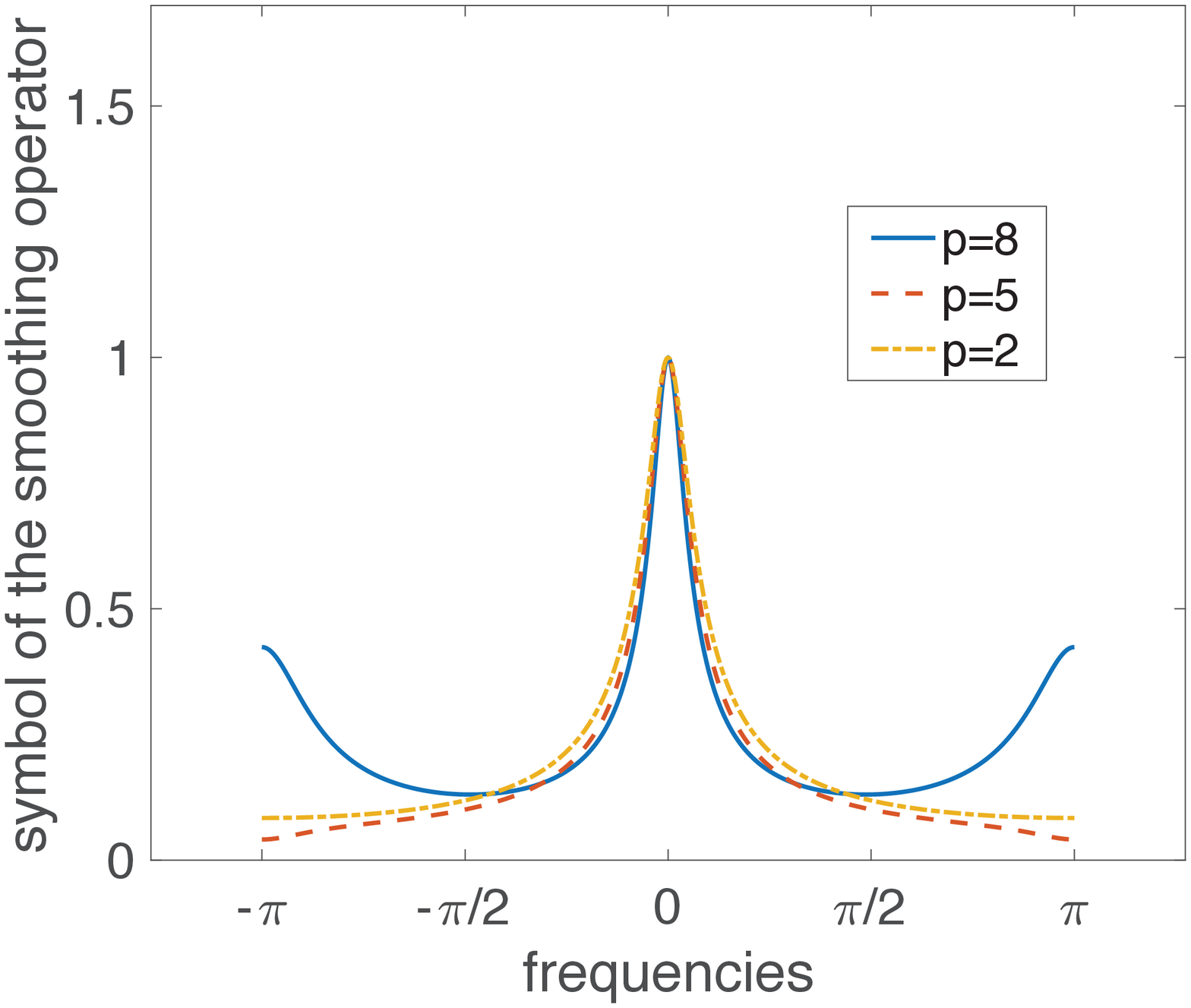}
\\
(a) Three-point Schwarz & (b) Five-point Schwarz
\end{tabular}
\label{symbol_Vankas}
\caption{Symbol of the smoothing operator corresponding to (a) the three-point and (b) the five-point multiplicative Schwarz smoothers, for three different spline degrees $p = 2, 5, 8$.}
\end{center}
\end{figure}

\begin{figure}[htb]
\begin{center}
\begin{tabular}{ccc}
\includegraphics[scale = 0.2]{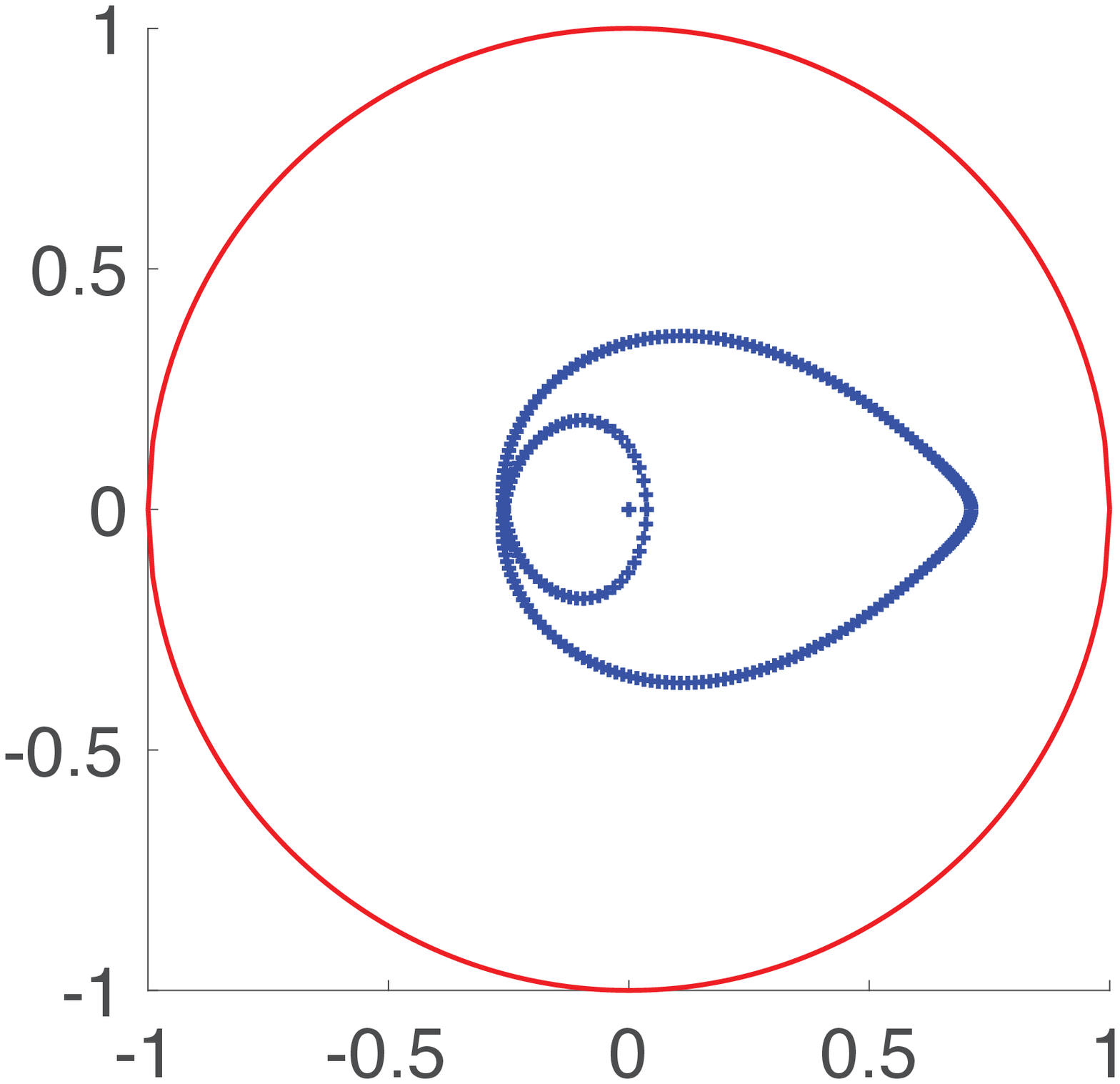}
&
\includegraphics[scale = 0.2]{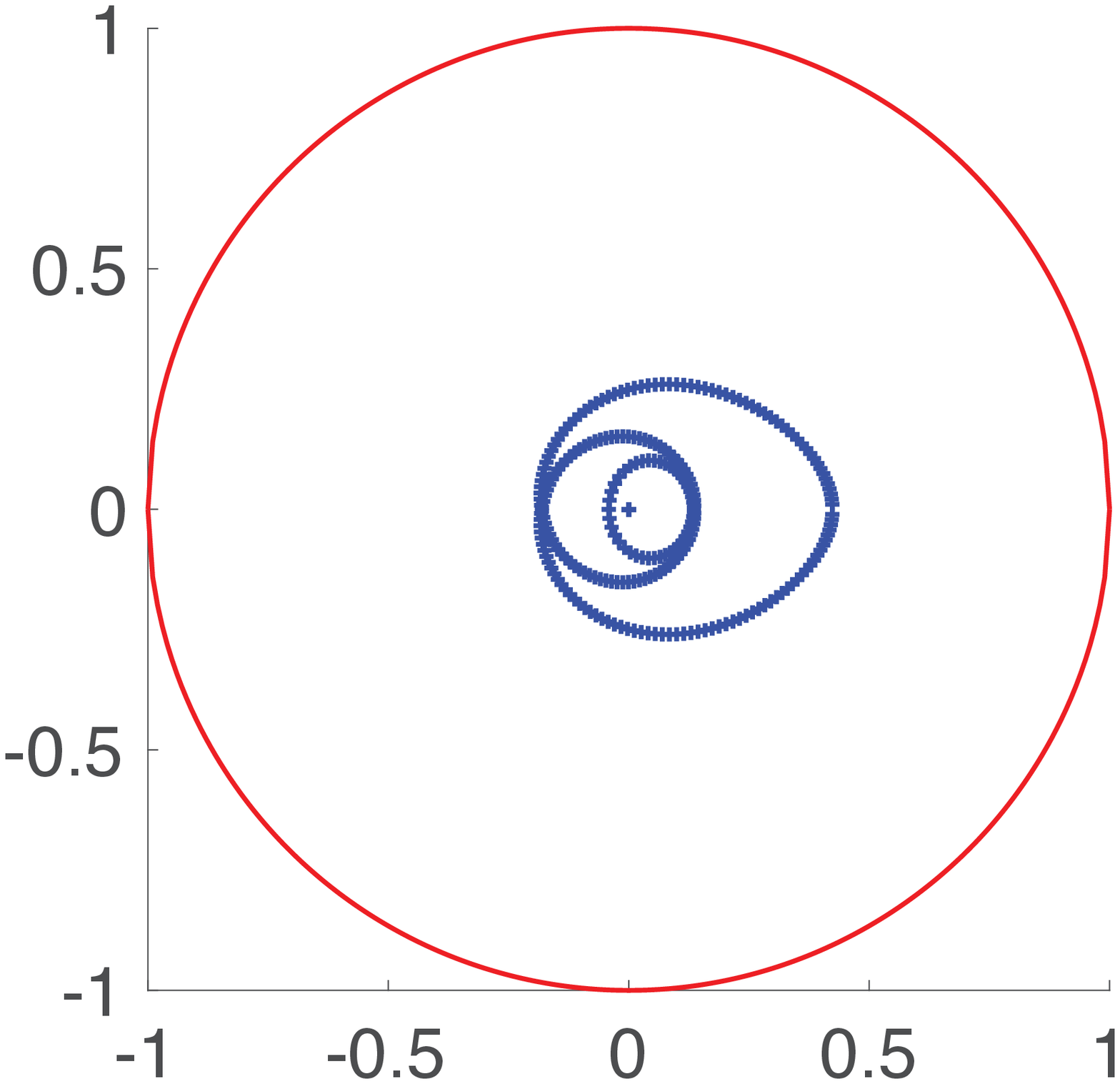}
&
\includegraphics[scale = 0.2]{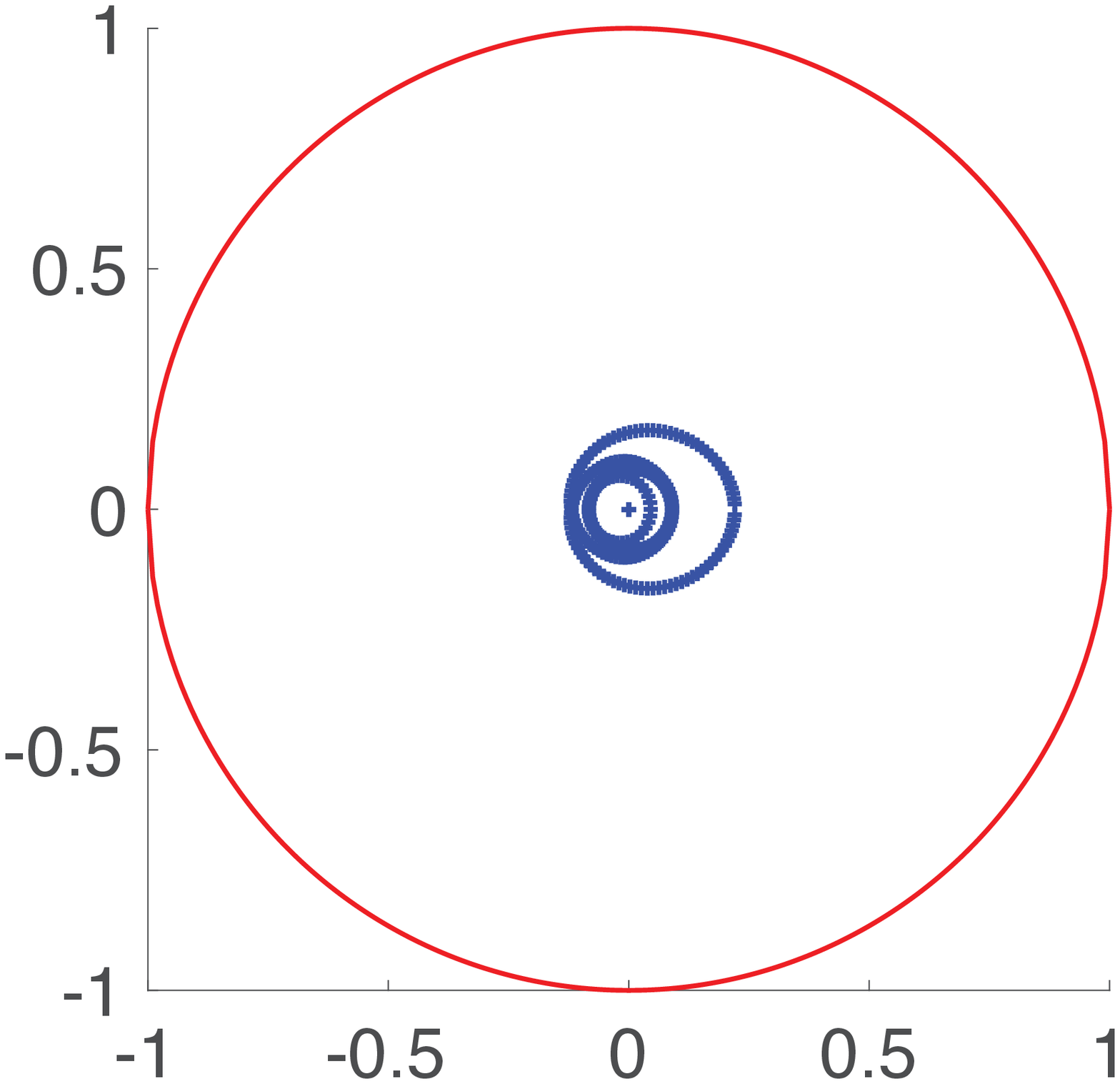}
\\
(a) Three-point Schwarz & (b) Five-point Schwarz & (c) Seven-point Schwarz
\end{tabular}
\caption{Distribution of the eigenvalues of the two-grid method based on the three different overlapping multiplicative Schwarz smoothers considered, when applied to $p=8$.}
\label{tg_eig_Vankas}
\end{center}
\end{figure}

This analysis can be generalized for any isogeometric discretization with spline degree $p$. The corresponding stencil has $(2p+1)$ elements and it has the following form:
$A_{h,p} =  \left[ a_p, a_{p-1},\ldots, a_1, a_0,a_1,\ldots,a_{p-1},a_p\right].$
Then, the matrices $P$ and $Q$  for the three-point multiplicative Schwarz smoother are given as follows:
$$P = \left(
\begin{array}{ccc}
a_2e^{\imath \theta} & a_1 & \sum_{j=0}^{p}a_je^{-\imath (j+1)\theta} \\ [1ex]
a_1e^{\imath \theta} & a_0 & \sum_{j=1}^{p}a_je^{-\imath j\theta} \\ [1ex]
a_0e^{\imath \theta} & a_1 & \sum_{j=2}^{p}a_je^{-\imath (j-1)\theta}\\
\end{array}
\right), \quad 
Q = \left(
\begin{array}{c}
-\sum_{j=3}^p a_j e^{\imath (j-1)\theta} \\ [1ex]
-\sum_{j=2}^p a_j e^{\imath j\theta} \\ [1ex]
-\sum_{j=1}^p a_j e^{\imath (j+1)\theta}
\end{array}
\right).
$$
It is even more interesting to analyze multiplicative Schwarz iterations with larger blocks, since our strategy will be to increase the size of the blocks as the spline degree $p$ grows up. LFA for these relaxations can be done in the same way that for the three-point multiplicative Schwarz smoother, but with heavier computations. Matrices $P$ and $Q$ for an arbitrary $n-$point multiplicative Schwarz iteration are $(n\times n)-$ and $(n\times 1)-$matrices respectively given by
$$P = \overline{P} D_P, \quad  \hbox{and} \quad 
Q = \left( 
\begin{array}{c}
-\sum_{j=n}^{p}a_j e^{\imath \left(j-\frac{n-1}{2} \right)\theta} \\ 
-\sum_{j=n-1}^{p}a_j e^{\imath \left(j-\frac{n-1}{2}+1 \right)\theta} \\
\vdots \\
-\sum_{j=\frac{n-1}{2}+1}^{p}a_j e^{\imath j \theta} \\
-\sum_{j=\frac{n-1}{2}}^{p}a_j e^{\imath \left(j+1\right)\theta} \\
\vdots \\
-\sum_{j=1}^{p}a_j e^{\imath \left(j+\frac{n-1}{2} \right)\theta}
\end{array}
\right),$$
where matrix $\overline{P}$ is 
$$
\left(\!\!
\begin{array}{ccccccccc} 
a_{n-1} & a_{n-2} & a_{n-3} & \cdots & a_{\frac{n-1}{2}} & \cdots & a_2 & a_1 & \sum_{j=0}^{p}a_j e^{-\imath \left(j+\frac{n-1}{2} \right)\theta}\\
a_{n-2} & a_{n-3} & \cdots & & \vdots & \dotso & a_1 & a_0 & \sum_{j=1}^{p}a_j e^{-\imath \left(j+\frac{n-1}{2}-1 \right)\theta}\\
a_{n-3} & \cdots &  & & \vdots & \cdots & a_0 & a_1 &  \sum_{j=2}^{p}a_j e^{-\imath \left(j+\frac{n-1}{2}-2 \right)\theta}\\
\vdots &  & & & \vdots & & & & \vdots \\
\vdots &  &  &  & a_0 &  &  &  & \sum_{j=\frac{n-1}{2}}^{p}a_j e^{-\imath j \theta}\\
\vdots & &  & &  \vdots & &. & & \vdots \\ 
a_2 & a_1 & a_0 & \cdots & \ldots &  & \cdots  & a_{n-4} & \sum_{j=n-3}^{p}a_j e^{-\imath \left(j-\frac{n-1}{2}+2 \right)\theta}\\
a_1 & a_0 & a_1 & \cdots  & \ldots & \cdots  & a_{n-4} & a_{n-3} & \sum_{j=n-2}^{p}a_j e^{-\imath \left(j-\frac{n-1}{2}+1 \right)\theta}\\
a_0 & a_1 & \cdots & \cdots & a_{\frac{n-1}{2}} & \cdots &  a_{n-3} & a_{n-2} & \sum_{j=n-1}^{p}a_j e^{-\imath \left(j-\frac{n-1}{2} 
\! \right)\theta}\\
\end{array}
\!\! \right),$$ 
and $D_P$ is the following diagonal matrix
$$D_P = \diag \left\{ \!
\begin{array}{ccccccccc}
e^{\imath \frac{n-1}{2}\theta}, &  e^{\imath \left(\frac{n-1}{2}-1\right)\theta}, & \cdots , & e^{\imath \theta},  & 1, & e^{-\imath \theta}, & \cdots , &  e^{-\imath \left(\frac{n-1}{2}-1\right)\theta},  & 1 
\end{array}
\!\right\}.$$
Thus, the presented LFA can be used for the choice of an adequate multiplicative Schwarz smoother for each spline degree. In particular, as we did previously for Gauss-Seidel relaxation, we can study if the class of multiplicative Schwarz smoothers are able to eliminate the high-frequency components of the error for IGA with a large spline degree. In Figure \ref{symbol_Vankas} we show the symbol of the smoothing operator of the three-point and five-point Schwarz relaxations  for each considered $p$. For $p=2$, we can see that the three-point approach is enough to adequately remove the high components. This approach, however, provides worse results for $p=8$. If, alternatively, we choose the five-point multiplicative Schwarz smoother, we can see in the picture that the high frequencies are removed for $p=2$ and $p=5$ in a very efficient way, but for $p=8$ it is not yet very satisfactory, and it would be recommended to increase the size of the blocks. 
We also analyze the eigenvalues of the two-grid operator based on the overlapping multiplicative Schwarz smoothers. In Figure \ref{tg_eig_Vankas} we can see for a fixed spline degree, $p=8$, that the eigenvalues become smaller when the size of the blocks within the Schwarz iteration gets larger. 
For an arbitrary spline degree $p$, this local Fourier analysis gives us the asymptotic convergence rates of the multigrid method based on the different multiplicative Schwarz smoothers. Thus, this information, together with the computational cost of each relaxation procedure can be used to justify the proposed strategy of varying the block size within the smoother for different values of $p$. 


\subsection{Local Fourier analysis results}
In order to demonstrate the utility of the described local Fourier analysis, next we show the good agreement between the convergence factors predicted by the LFA and the asymptotic convergence factors experimentally obtained by the multigrid method. With this purpose, in Table \ref{table_LFA_1D} we provide the smoothing ($\rho_{1g}$), two-grid ($\rho_{2g}$) and three-grid ($\rho^V_{3g}$) convergence factors obtained from LFA, considering one smoothing step of a Gauss-Seidel relaxation, together with the asymptotic convergence factors provided by the $W(1,0)-$ and $V(1,0)-$cycle multigrid codes ($\rho^W_h$ and $\rho^V_h$, respectively). We can observe a perfect match between the experimental factors and those predicted by LFA for both W- and V-cycles. 
\begin{table}[htb]
\begin{center}
\begin{tabular}{cccccc }
\cline{2-6}
&  & & Gauss-Seidel \\
\cline{2-6}
& \quad$\rho_{1g}$\quad  & \quad$\rho_{2g}$\quad & \quad$\rho^W_h$\quad & \quad$\rho^V_{3g}$\quad & \quad$\rho^V_h$\quad \\
\hline
$p = 2$ & 0.31 & 0.19 & 0.19 & 0.19 & 0.19 \\
$p = 3$ & 0.26 & 0.22 & 0.22 & 0.22 & 0.22 \\
$p = 4$ & 0.38 & 0.38 & 0.38 & 0.38 & 0.38 \\
$p = 5$ & 0.62 & 0.62 & 0.62 & 0.62 & 0.62 \\
$p = 6$ & 0.79 & 0.79 & 0.80 & 0.79 & 0.80 \\
$p = 7$ & 0.89 & 0.89 & 0.90 & 0.89 & 0.90 \\
$p = 8$ & 0.99 & 0.99 & 0.96 & 0.99 & 0.96 \\
\hline
\end{tabular}
\caption{Smoothing ($\rho_{1g}$), two-grid ($\rho_{2g}$) and three-grid ($\rho^V_{3g}$) convergence factors predicted by LFA together with the asymptotic convergence factors provided by the W(1,0) and V(1,0) cycle multigrid codes ($\rho^W_h$ and $\rho^V_h$, respectively), for different spline degrees $p$.}
\label{table_LFA_1D}
\end{center}
\end{table}

We also perform the analysis presented in Section \ref{sec:4_schwarz} for the multiplicative Schwarz smoothers. In Table \ref{table_LFA_Schwarz_1D}, we show the smoothing and three-grid convergence factors predicted by the analysis and the experimentally obtained asymptotic convergence factors of a $V-$cycle multigrid with one pre- and no post-smoothing steps. These results are shown for the three considered multiplicative Schwarz smoothers (three-, five- and seven-point approaches) and for different spline degrees from $p=2$ to $p=8$. We consider $V-$cycles because we have seen that their convergence rates are as those provided by $W-$cycles.
We observe in the table that the smoothing ability of the proposed three-point multiplicative Schwarz relaxation deteriorates when $p$ becomes larger. 
This affects the three-grid convergence factor that also gets worse, showing the necessity of considering multiplicative Schwarz smoothers coupling more than three points. 
\begin{table}[htb]
\begin{center}
\begin{tabular}{cccccccccc}
\cline{2-10}
& \multicolumn{3}{c}{3p Schwarz} & \multicolumn{3}{c}{5p Schwarz} & \multicolumn{3}{c}{7p Schwarz}\\
\cline{2-10}
& $\rho_{1g}$  & $\rho^V_{3g}$ & $\rho^V_h$ & $\rho_{1g}$  & $\rho^V_{3g}$ & $\rho^V_h$ & $\rho_{1g}$  & $\rho^V_{3g}$ & $\rho^V_h$ \\
\hline
$p = 2$ & 0.176 & 0.127 & 0.127 & 0.119 & 0.088 & 0.087 & 0.089 & 0.065 & 0.065 \\
$p = 3$ & 0.156 & 0.114 & 0.113 & 0.112 & 0.086 & 0.086 & 0.086 & 0.066 & 0.066 \\
$p = 4$ & 0.146 & 0.127 & 0.127 & 0.104 & 0.084 & 0.084 & 0.082 & 0.067 & 0.067 \\
$p = 5$ & 0.209 & 0.209 & 0.211 & 0.101 & 0.095 & 0.095 & 0.078 & 0.069 & 0.069 \\
$p = 6$ & 0.389 & 0.389 & 0.389 & 0.147 & 0.147 & 0.147 & 0.077 & 0.077 & 0.077 \\
$p = 7$ & 0.564 & 0.564 & 0.564 & 0.279 & 0.279 & 0.276 & 0.119 & 0.119 & 0.121 \\
$p = 8$ & 0.712 & 0.712 & 0.712 & 0.424 & 0.424 & 0.426 & 0.221 & 0.221 & 0.224 \\
\hline
\end{tabular}
\caption{Smoothing ($\rho_{1g}$) and three-grid ($\rho^V_{3g}$) convergence factors predicted by LFA together with the asymptotic convergence factors provided by the V(1,0)-cycle multigrid code ($\rho^V_h$), for different spline degrees $p$.}
\label{table_LFA_Schwarz_1D}
\end{center}
\end{table}
We have demonstrated that the local Fourier analysis is a very useful tool to obtain information about the performance of multigrid for IGA, since it predicts very accurately the convergence of the method. This is very interesting since for a fixed spline degree $p$, we can choose the appropriate number of points in the blocks to construct an efficient multiplicative Schwarz smoother. 
\subsubsection{Smoothing strategy}\label{sec:4_strategy}
\begin{figure}[htb]
\begin{center}
\includegraphics[scale = 0.32]{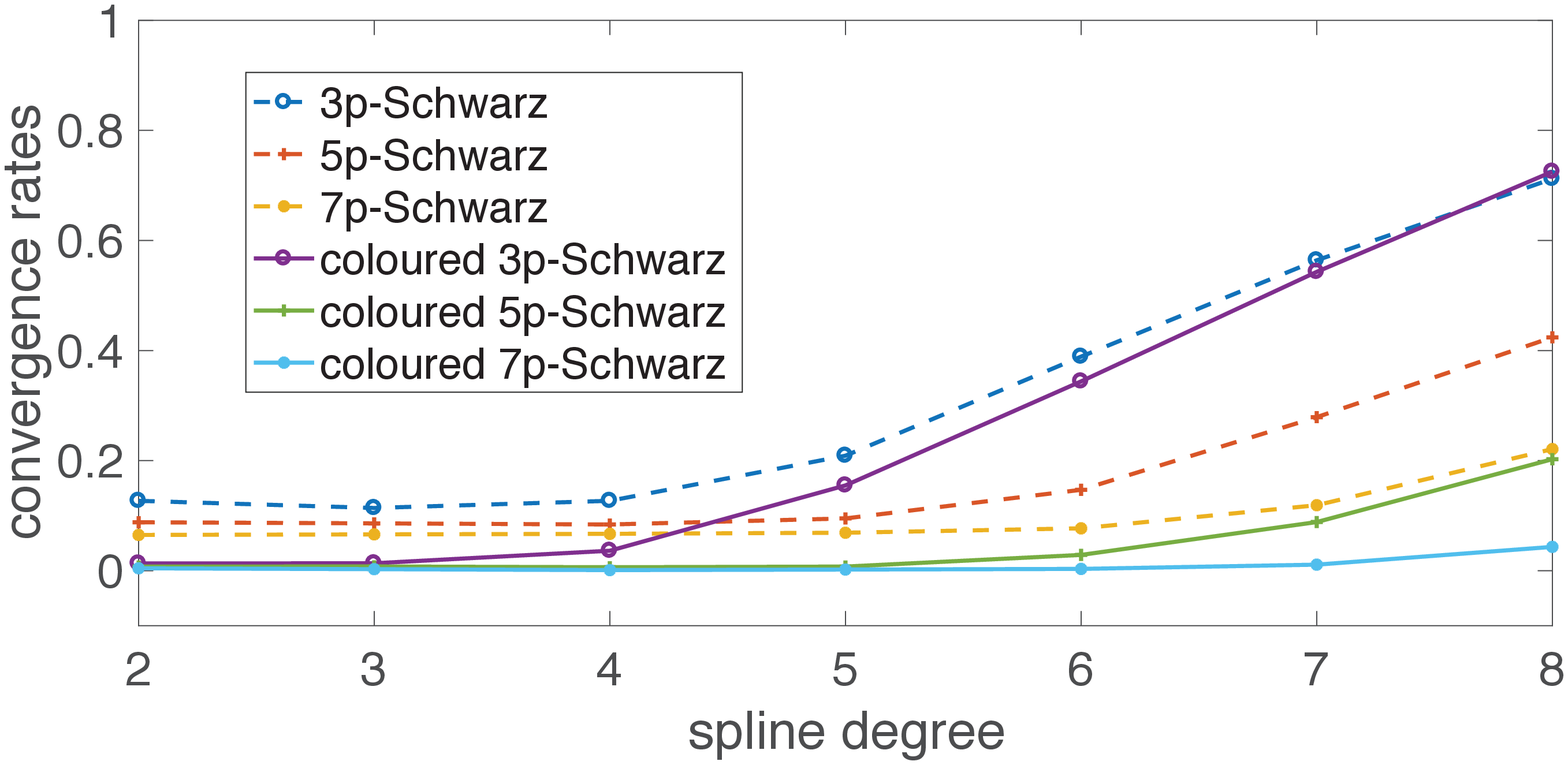}
\label{schwarz_rates}
\caption{Convergence rates of the multigrid method based on overlapping multiplicative Schwarz smoothers for the three considered Schwarz relaxations and their coloured counterparts, and for spline degrees $p=2,\ldots, 8$.}
\end{center}
\end{figure}

We show in Figure \ref{schwarz_rates}, by using dotted lines, the convergence rates provided by the proposed multigrid method based on the three-, five-, and seven-point multiplicative Schwarz smoothers, for $p = 2,\ldots, 8$. We can observe that $p = 4$ and $p=6$ seem to be a point of change in the behavior of the three- and five-point smoothers, so, this can be taken as an advice for the choice of the smoother. 

Notice also that the different blocks within the Schwarz smoother can also be visited in different orderings, for instance, they can also be treated with some patterning scheme, yielding a multicolored version of these relaxation schemes. This type of approaches usually provides better convergence. In our case, we have done a comparison of the convergence rates of the considered multiplicative Schwarz smoothers and their coloured counterparts. In particular, we have considered a three-colour version of these relaxation procedures, and we have seen that the number of iterations necessary to achieve the desired stopping criterium are reduced almost in half when the coloured versions are considered. For this reason, it seems interesting to consider these latter approaches. Thus, in Figure \ref{schwarz_rates}, by using continuous lines, we also display the asymptotic convergence factors provided by a multigrid based on $V(1,0)-$cycles and the coloured versions of the considered three-, five- and seven-point multiplicative Schwarz smoothers. We can observe that the qualitative behavior of the method is as in the case for the non-coloured smoothers, but the convergence rates are much better. Taking into account the convergence rates in Figure \ref{schwarz_rates} and some CPU times, we can fix that the most efficient strategy is the coloured three-point multiplicative Schwarz smoother for $p= 2, 3$, its five-point extension for $p = 4, 5$, and finally the Schwarz relaxation with seven-point size blocks for $p = 6, 7, 8$. 
For the two-dimensional case, we extend this strategy directly and we have obtained again very satisfactory results, as we will see in the numerical experiments presented in the next section.


\section{Numerical results}\label{sec:5}
\setcounter{section}{5}
In this section, we consider three different numerical experiments to illustrate the 
robustness and efficiency of the proposed multigrid method based on multiplicative Schwarz smoothers. First, we will solve both one- and two-dimensional problems on the parametric space. After that, we will deal with a two-dimensional problem on a nontrivial computational domain. For the first two numerical experiments we choose B-splines as basis functions, whereas for the third numerical experiment NURBS are used because these functions can exactly describe the geometry for the considered domain. 

Since the results obtained by the three-grid Fourier analysis in Section \ref{sec:4} demonstrated that V-cycles provide similar convergence as W-cycles, the former are chosen because they have a lower computational cost and therefore a more efficient multigrid method is obtained. Moreover, only one smoothing step is considered in all the numerical experiments, in particular V(1,0)-cycles are used. In all the cases the initial guess is taken as a random vector, and the stopping criterion for the multigrid solver is set to reduce the initial residual by a factor of $10^{-8}$. All numerical computations were carried out using MATLAB on a MacBook Pro with a Core i5 2.7 GHz and 8 GB RAM, running OS X 10.10 (Yosemite).

\subsection{One-dimensional example on the parametric domain}
In the first numerical experiment we consider the following two-point boundary value problem 
$$
\left \{
\begin{array}{l}
-u''(x) = \pi^2 \sin(\pi x), \quad x \in (0,1), \nonumber \\
u(0) = u(1) = 0. \nonumber
\end{array}
\right.
$$
We discretize this problem by using an equidistant knot span and maximum continuity splines for different degrees ranging from $p=2$ until $p=8$. For solving 
the linear system we apply the proposed $V(1,0)$ multigrid cycle based on the coloured multiplicative Schwarz 
smoother. The size of the blocks in the relaxation depends on the spline degree and it has been chosen following the results provided by the local Fourier analysis. As commented in Section \ref{sec:4}, we consider blocks composed of three points for the cases $p=2, 3$, five points for $p=4, 5$ and seven points for $p=6,7, 8$. 
Table \ref{tab:Sch_1D} shows the number of iterations needed to reduce the initial residual by a factor of $10^{-8}$ for several mesh sizes and different spline degrees  $p = 2, \ldots, 8$. We also report the computational time (in seconds) needed for solving the linear system by the proposed multigrid method. 

 
  \begin{table}[htbp]
 	\label{tab:Sch_1D}
 	\centering
 	\begin{tabular}{c|c|c|c|c|c|c|c|}  
	\cline{2-8}
	& \multicolumn{2}{|c|}{Color 3p Schwarz} & \multicolumn{2}{|c|}{Color 5p Schwarz} & \multicolumn{3}{|c|}{Color 7p Schwarz}\\
	 \cline{2-8} 
 	& $p=2$ & $p=3$ & $p=4$ & $p=5$ & $p=6$ & $p=7$ & $p=8$  \\ \hline 
	\multicolumn{1}{|c|}{$h^{-1} $}  & it \; cpu & it \; cpu &  it \; cpu & it \; cpu & it \; cpu & it \; cpu & it \; cpu \\
	\hline
         \multicolumn{1}{|c|}{512} & $ 5 \quad 0.63$ & $ 5 \quad 0.75$ & $ 4 \quad 0.76$ & $ 4 \quad 0.80$
	 & $ 4 \quad 0.95$ & $ 4 \quad 1.07$ & $ 5 \quad 1.42$ \\ 
	\multicolumn{1}{|c|}{1024} & $ 5 \quad 0.93$ & $ 5 \quad 1.04$ & $ 4 \quad 1.11$ & $ 4 \quad 1.16$
	 & $ 4 \quad 1.39$ & $ 4 \quad 1.47$ & $ 5 \quad 1.90$ \\ 
	\multicolumn{1}{|c|}{2048} & $ 5 \quad 1.63$ & $ 5 \quad 1.75$ & $ 4 \quad 1.88$ & $ 4 \quad 1.94$
	 & $ 4 \quad 2.34$ & $ 4 \quad 2.49$ & $ 5 \quad 3.20$ \\ 
	\multicolumn{1}{|c|}{4096} & $ 5 \quad 3.00$ & $ 5 \quad 3.18$ & $ 4 \quad 3.43$ & $ 4 \quad 3.51$
	 & $ 4 \quad 4.31$ & $ 4 \quad 4.47$ & $ 5 \quad 5.72$ \\ 
	\multicolumn{1}{|c|}{8192} & $ 5 \quad 5.83$ & $ 5 \quad 6.17$ & $ 4 \quad 6.60$ & $ 4 \quad 6.79$
	 & $ 4 \quad 8.24$ & $ 4 \quad 8.60$ & $ 5 \quad 10.77$ 
 \\ \hline
 	\end{tabular}
 	\caption{One-dimensional test problem. Number of $V(1,0)$ multigrid iterations (it) and computational time (cpu) necessary to reduce the initial residual in a factor of $10^{-8},$ for different grid-sizes $h$ and for different values of the spline degree $p$, using the most appropriate coloured multiplicative Schwarz smoother for each $p$.}	
 \end{table}
 
First of all, we observe that our approach shows an excellent convergence with respect to the mesh refinement, for all the polynomial orders. But overall, we would like to remark the robustness of the method with respect to the spline degree. In all cases, only four or five $V(1,0)-$cycles are needed to reach the stopping criterium, independently of the mesh size and the polynomial order. Regarding the computational cost, we observe a very good scaling, not only with respect to $h$ but also with respect to $p$. 

\subsection{Two-dimensional example on the parametric domain}
We now test the performance of the geometric multigrid solver based on overlapping multiplicative Schwarz smoothers on a two-dimensional problem defined on the parametric domain $\Omega = (0,1)^2$. More concretely, we consider the problem,
$$
\left \{
\begin{array}{l}
- \Delta u = 2 \pi^2 \sin(\pi x) \sin(\pi y), \quad (x,y) \in \Omega = (0,1)^2, \nonumber \\
u(x,y) = 0, \quad (x,y) \; {\mbox on} \; \partial \Omega \nonumber
\end{array}
\right.
$$
Since the problem is solved in the parametric domain, we choose the basis functions as B-splines.
Similar to the one-dimensional case, the performance of the geometric multigrid solver based on classical smoothers like Gauss-Seidel iteration, is very poor. For example, for a mesh $32\times 32$ and a moderate large spline degree $p=4$, the convergence factor of this standard multigrid solver is already close to one, requiring more than $700$ iterations to reach the convergence.
This is due to the presence of many small eigenvalues associated with oscillatory components of the error. 

We now apply the multigrid method based on the proposed coloured overlapping multiplicative Schwarz smoothers. The size of the blocks depend on the spline degree. In particular, based on the Fourier analysis results in Section \ref{sec:4}, we choose blocks of size $3 \times 3$ for the cases $p=2, 3$, blocks of size $5\times 5$ for the cases $p=4, 5$ and  blocks of size $7\times 7$ for splines degree $p=6, 7, 8$. Again, we use $V-$cycles with only one pre-smoothing step and no post-smoothing steps. The number of iterations needed to reduce the initial residual by a factor of $10^{-8}$ for several mesh sizes and different spline degrees  $p = 2, \ldots, 8$, together with the computational time (in seconds) are given in Table \ref{tab:Sch_2D}.  We observe that the iteration numbers are robust with respect to both the size of the grid $h$ and the spline degree $p$. Moreover, we see that the number of iterations is similar to those reported in the one-dimensional case. We emphasize that only four or five $V(1,0)-$cycles are needed to reach the stopping criterium, independently of $h$ and $p$. Just like in the one-dimensional case, we can conclude that the multigrid method based on an appropriate multiplicative Schwarz smoother provides an efficient and robust solver for B-spline isogeometric discretizations.
  \begin{table}[htbp]
 	\label{tab:Sch_2D}
 	\centering
 	\begin{tabular}{c|c|c|c|c|c|c|c|}  
	\cline{2-8}
	& \multicolumn{2}{|c|}{Color 3p Schwarz} & \multicolumn{2}{|c|}{Color 5p Schwarz} & \multicolumn{3}{|c|}{Color 7p Schwarz}\\
	 \cline{2-8} 
 	& $p=2$ & $p=3$ & $p=4$ & $p=5$ & $p=6$ & $p=7$ & $p=8$  \\ \hline 
	\multicolumn{1}{|l|}{Grid}  & it \; cpu & it \; cpu &  it \; cpu & it \; cpu & it \; cpu & it \; cpu & it \; cpu \\
	\hline
         \multicolumn{1}{|l|}{$32^2$} & $ 4 \quad 0.55$ & $ 4 \quad 0.57$ & $ 3 \quad 0.82$ & $ 4 \quad 1.13$
	 & $ 3 \quad 1.43$ & $ 3 \quad 1.66$ & $ 4 \quad 2.26$ \\ 
	\multicolumn{1}{|l|}{$64^2$} & $ 4 \quad 1.25$ & $ 4 \quad 1.31$ & $ 3 \quad 2.41$ & $ 4 \quad 3.53$
	 & $ 3 \quad 5.27$ & $ 3 \quad 5.36$ & $ 5 \quad 9.59$ \\ 
	\multicolumn{1}{|l|}{$128^2$} & $ 4 \quad 5.30$ & $ 4 \quad 5.59$ & $ 3 \quad 9.20$ & $ 4 \quad 13.45$
	 & $ 3 \quad 18.01$ & $ 3 \quad 19.36$ & $ 5 \quad 36.71$ \\ 
	\multicolumn{1}{|l|}{$256^2$} & $ 4\quad 28.41$ & $4 \quad 28.71$ & $ 3 \quad 41.27$ & $4 \quad 54.42$
	 & $3  \quad 78.88$ & $ 3 \quad 79.57$ & $ 5 \quad 141.16$ 
 \\ \hline
 	\end{tabular}
 	\caption{Two-dimensional test problem. Number of $V(1,0)$ multigrid iterations (it) and computational time (cpu) necessary to reduce the initial residual in a factor of $10^{-8},$ for different grid-sizes and for different values of the spline degree $p$, using the most appropriate coloured multiplicative Schwarz smoother for each $p$.}	
 \end{table}

\subsection{Quarter annulus}

\begin{figure}[htb]
\begin{center}
\includegraphics[scale = 0.7]{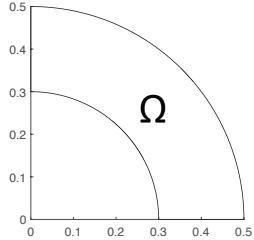}
\label{annulus_domain}
\caption{Computational domain for the Poisson problem on a quarter of an annulus.}
\end{center}
\end{figure}

The last experiment demonstrates the efficiency and robustness of the proposed multigrid method to deal with a nontrivial geometry. We take as domain the quarter of an annulus, 
$$
\Omega = \{(x,y) \in \mathbb{R}^2 \, | \, r^2 \leq x^2 + y^2 \leq R^2, x, y \geq 0\},
$$
which is sketched in Figure \ref{annulus_domain}. We consider the solution of the Poisson problem in such domain with homogeneous Dirichlet boundary conditions
$$
\left \{
\begin{array}{l}
- \Delta u = f(x,y), \quad (x,y) \in \Omega, \nonumber \\
u(x,y) = 0, \quad (x,y) \; {\mbox on} \; \partial \Omega, \nonumber
\end{array}
\right.
$$
where $f(x,y)$ is such that the exact solution is 
$$
u(x,y) = \sin(\pi x) \sin(\pi y) (x^2+y^2-r^2)(x^2+y^2-R^2).
$$
The geometry of the computational domain is described exactly by quadratic $C^1$ NURBS, and we
choose in our experiments $r=0.3$ and $R=0.5$. To discretize this problem, we use NURBS of
degree $p=2,\ldots,8$ with maximal smoothness. We solve the corresponding linear systems using $V(1,0)-$ cycles based on the coloured multiplicative Schwarz smoother described in Section \ref{sec:4_strategy}. The size of the blocks of the Schwarz smoother depends on $p$, and they are chosen following the same strategy  than in the previous experiment. The numbers of iterations needed to reach the stopping criterium, for various degrees $p$ and for different mesh sizes, are reported in Table \ref{tab:Sch_2D_A}. We observe an excellent performance, obtaining similar results to those of the example with the parametric domain. We can conclude that the proposed solver is robust with respect to this geometry transformation, having a great potential for solving problems in more complicated multi-patch geometries, which is subject of future research.
 
  \begin{table}[htbp]
 	\label{tab:Sch_2D_A}
 	\centering
 	\begin{tabular}{c|c|c|c|c|c|c|c|}  
	\cline{2-8}
	& \multicolumn{2}{|c|}{Color 3p Schwarz} & \multicolumn{2}{|c|}{Color 5p Schwarz} & \multicolumn{3}{|c|}{Color 7p Schwarz}\\
	 \cline{2-8} 
 	& $p=2$ & $p=3$ & $p=4$ & $p=5$ & $p=6$ & $p=7$ & $p=8$  \\ \hline 
         \multicolumn{1}{|l|}{$32^2$} & $ 4 $ & $ 4 $ & $ 3 $ & $ 4 $
	 & $ 3 $ & $ 3 $ & $ 4 $ \\ 
	\multicolumn{1}{|l|}{$64^2$} & $ 4 $ & $ 4 $ & $ 3$ & $ 4 $
	 & $ 3 $ & $ 3 $ & $ 5 $ \\ 
	\multicolumn{1}{|l|}{$128^2$} & $ 4 $ & $ 4 $ & $ 3 $ & $ 4 $
	 & $ 3 $ & $ 3$ & $ 5 $ \\ 
	\multicolumn{1}{|l|}{$256^2$} & $ 4$ & $4 $ & $ 3 $ & $4$
	 & $3 $ & $ 3 $ &  $5 $ 
 \\ \hline
 	\end{tabular}
 	\caption{Quarter annulus problem. Number of $V(1,0)$ multigrid iterations (it) necessary to reduce the initial residual in a factor of $10^{-8},$ for different grid-sizes and for different values of the spline degree $p$, using the most appropriate coloured multiplicative Schwarz smoother for each $p$.}	
 \end{table}



\section{Conclusions}\label{sec:6}
\setcounter{section}{6}

An efficient and robust geometric multigrid method, based on overlapping multiplicative Schwarz smoothers, is proposed for IGA. The robustness of the algorithm with respect to the spline degree is demonstrated through numerical experiments and a local Fourier analysis. This analysis is carried out in 1D for any spline degree $p$ and an arbitrary size of the blocks to solve in the relaxation procedure. The key point to achieve a robust algorithm is the choice of larger blocks within the Schwarz smoother when the spline degree grows up. Moreover, to improve the convergence rates provided by the multigrid method, we considered a coloured version of the multiplicative Schwarz smoothers.  A simple multigrid V-cycle based on this type of smoothers with only one smoothing step results to be a robust and efficient algorithm for isogeometric discretizations.


%

\bibliographystyle{plain}
\bibliography{references}
\end{document}